\documentclass[11pt]{amsart}

\usepackage{epsf,amssymb}

\newtheorem{thm}{Theorem}[section]
\newtheorem{lem}[thm]{Lemma}
\newtheorem{cor}[thm]{Corollary}
\newtheorem{prop}[thm]{Proposition}
\newtheorem{rem}[thm]{Remark}
\newtheorem{strategy}[thm]{Strategy}

\newtheorem*{ttm}{Theorem}
\newtheorem{quest}{Question}

\newtheorem*{claim}{Claim}

\newcommand{\bbr}{\begin{rem}\em} 
\newcommand{\eer}{\end{rem}}
\newcommand{\bdry}{\partial}
\newcommand{\m}{mult}

\newcommand{\s}{\vskip.1in}

\newcommand{\be}{\begin{enumerate}}
\newcommand{\ee}{\end{enumerate}}
\renewcommand{\backslash}{\setminus}

\def\tb{\operatorname{tb}}
\def\C{\hbox{$\mathbb C$} }
\def\Z{\hbox{$\mathbb Z$} }
\def\Q{\hbox{$\mathbb Q$} }
\def\R{\hbox{$\mathbb R$} }

\def\co{\colon\thinspace}

\def\dfn#1{{\em #1}}

\begin{document}

\title{Knots and Contact Geometry}
\date{June 15, 2000}

\author{John B. Etnyre}
\address{Stanford University, Stanford, CA 94305}
\email{etnyre@math.stanford.edu}
\urladdr{http://math.stanford.edu/\char126 etnyre}

\author{Ko Honda}
\address{University of Georgia, Athens, GA 30602}
\email{honda@math.uga.edu}
\urladdr{http://www.math.uga.edu/\char126 honda}


\begin{abstract}
	We classify Legendrian torus knots and Legendrian figure eight knots in the tight contact structure on $S^3$
	up to Legendrian isotopy. As a corollary to this we also obtain the classification of transversal 
	torus knots and transversal figure eight knots up to transversal isotopy.
\end{abstract}

\maketitle


\section{Introduction}

There have recently been several breakthroughs concerning the classification of
tight contact structures on 3-manifolds \cite{Giroux99, Honda1, Honda2, Honda3, EH}. The main ingredient
in all these advances is the theory of {\em{convex surfaces}}, due to Giroux \cite{gi:convex}.
Convexity enables us
to understand neighborhoods of surfaces very well and, more surprisingly, relate the
characteristic foliation on one surface to the characteristic foliation on another
surface that meets it along a Legendrian curve. Moreover, it reduces the study of
characteristic foliations on a surface --- often quite delicate ---  to the study
of multi curves on a surface.  These {\it dividing curves} turn out to be more flexible and robust.

In this paper the authors introduce the tools from convex surface theory to the study 
of Legendrian and transversal knots in the standard tight contact structure on $S^3.$
Legendrian and transversal knots have played an important role in distinguishing tight
contact structures \cite{k,LM} and detecting overtwisted contact structures \cite{be,EH}.
Moreover, important topological properties can be detected using Legendrian knots. Specifically,
Rudolph \cite{r} has shown how to use invariants of Legendrian knots in a knot type to find 
obstructions to slicing this knot type. This has been generalized by Lisca and Mati\'c \cite{LM2} 
and Kronheimer and Mrowka \cite{KM} to find bounds on the slice genus of a knot. 

Recently there has been some progress in the classification of transversal knots. Eliashberg \cite{el:knots}
had previously shown that transversal unknots are determined by their self-linking number (the only classical invariant),
and, for a few years, this was the only knot type for which a transversal classification existed.
Then, in 1998,
Etnyre \cite{Etnyre99} classified positive transversal torus knots by showing their knot
type and self-linking number determine the transversal isotopy class.
More recently, Birman and Wrinkle \cite{bw}
extended Etnyre's results to transversal iterated torus knots, using a different approach  ---
namely, the study of braid foliations and the work of Menasco \cite{M} on iterated torus knots.
Until now, however, the only classification result for
Legendrian knots was the classification of Legendrian unknots by Eliashberg and Fraser \cite{ef}.
Eliashberg and Fraser had proved that, for Legendrian unknots,
the classical invariants (the Thurston-Bennequin invariant and the rotation
number) determine the Legendrian isotopy type. In this paper we prove:
\begin{ttm}
	Two oriented Legendrian torus knots are Legendrian isotopic if and only if 
	their Thurston-Bennequin invariants, rotation numbers and knot types agree.
\end{ttm}
\noindent
We also find the range of the classical invariants for Legendrian torus knots,
thus finishing the classification. In particular we show that, for a negative $(p,-q)$-torus knot  $K$ with $p>q>0$,
$$\tb(K)\leq -pq.$$
The classical Bennequin inequality only gives
$$\tb(K)\leq pq-p-q,$$
while the bounds discovered by Fuchs and Tabachnikov \cite{ft, Ta} give
$$\tb(K)\leq -pq$$
if $q$ is even but only
$$\tb(K)\leq -pq+p-q$$
when $q$ is odd. We thank Fuchs for informing us of the computations of these bounds in the dissertation \cite{ep}
of one of his students.  Thus we give the first class of knots --- the $(p,-q)$-torus knots with $q$ odd --- for
which all known bounds on the Thurston-Bennequin invariant are not sharp.

As a corollary of the above theorem, we obtain a special case of the Birman-Wrinkle-Menasco result:
\begin{ttm}
	Two transversal torus knots are transversally isotopic if and only if 
	their self-linking numbers and knots types agree.
\end{ttm}

We also prove:
\begin{ttm}
	Two oriented Legendrian figure eight knots are Legendrian isotopic if and only if
	their Thurston-Bennequin invariants and rotation numbers agree.
\end{ttm}
\noindent
We complete our classification by identifying the range of the Thurston-Bennequin invariant and rotation number for
Legendrian figure eight knots. We then obtain the classification of transversal figure eight knots as a corollary.
\begin{ttm}
	Two transversal figure eight knots are transversally isotopic if and only if
	their self-linking numbers agree.
\end{ttm}

\noindent
Continuing our line of inquiry, an open-ended question is:

\begin{quest}
Which transversal and Legendrian knots are determined by their classical invariants?
\end{quest}

It is known that, at least for Legendrian knots, the answer is not `all knot types'.
Using the powerful new invariants of contact homology,
Hofer and Eliashberg and (independently) Chekanov \cite{Chekanov} found Legendrian knots whose
classical invariants agree but are not Legendrian isotopic. It is interesting to note that the knot type 
Chekanov uses is the first knot type in the standard knot tables (e.g. \cite{Ro}) not covered by
one of the above theorems!  For transverse knots, the authors do not know whether there exist transversal
knots which are not distinguished by their classical invariants.


\section{Basic contact geometry}

Recall a \dfn{contact structure} on a 3-manifold $M$ is a maximally nonintegrable plane field
$\xi.$ Throughout this paper we assume our contact structures are transversely oriented and
thus can be globally given as the kernel of a 1-form $\alpha$ where $\alpha\wedge d\alpha\not=0.$
Moreover we always orient $M$ by $\alpha\wedge d\alpha.$

If $\Sigma$ is a surface in $M$ then $\xi\cap T\Sigma$ is a singular line field on $\Sigma$ 
and may be integrated to a singular foliation $\Sigma_\xi$ called the \dfn{characteristic 
foliation}. The singularities may be assumed to be either elliptic or hyperbolic (depending
on the local degree of the foliation) and
if $\Sigma$ is oriented then they also have a sign determined by the compatibility of the
orientations of $\xi$ and $T\Sigma$ at the singularities. There are many standard ways to manipulate
the characteristic foliation ---  for
details see \cite{a, gi:convex,ef}. We also recall that the characteristic foliation determines
a contact structure in a neighborhood of the surface.

A contact structure $\xi$ is called \dfn{tight} if there are no embedded disks $D$ 
with a limit cycle in their characteristic foliation. If $\xi$ is not tight, then it is called
\dfn{overtwisted}. The standard example of a tight contact structure is given by $\xi_0,$ the 
complex tangencies to $S^3\subset \C^2.$  This is the unique tight structure on $S^3$ \cite{el:twenty}.
The uniqueness is easily seen using Darboux's Theorem and the following theorem.
\begin{thm}[Eliashberg \cite{el:twenty}]\label{thm:uniqueB3}
	A tight contact structure on the 3-ball is uniquely determined (up to isotopy) by the 
	characteristic foliation on its boundary.
\end{thm}
In \cite{el:twenty} the group of contactomorphisms was also studied. Fix a point $p$ in $S^3$ and
let $\hbox{Diff}_0(S^3)$ be the
group of orientation-preserving diffeomorphisms of $S^3$ that fix the plane $\xi_0(p),$ and let
$\hbox{Diff}_{\xi_0}$ be the group of diffeomorphisms of $S^3$ that preserve $\xi_0.$ 
\begin{thm}[Eliashberg \cite{el:twenty}]\label{thm:htpyequiv}
	The natural inclusion of 
	$$\hbox{Diff}_{\xi_0}\hookrightarrow \hbox{Diff}_0(S^3)$$ is a
	weak homotopy equivalence.
\end{thm}

\subsection{Legendrian Knots}\label{legendriansection}

A  curve $\gamma$ in $M$ is called a \dfn{Legendrian curve} (or \dfn{Legendrian knot} if $\gamma$ is
closed) if it is everywhere tangent to the contact plane field $\xi.$  Suppose $\gamma$ is closed.
Our prime interest in this paper is the classification of a Legendrian knot $\gamma$ up to 
isotopy through Legendrian knots. This is equivalent to the classification of Legendrian knots
up to global contact isotopies.
\begin{lem}[see \cite{el:knots}]\label{extend}
	If $\psi_t\co S^1\to M$ is a Legendrian isotopy, then there
	is a contact isotopy $f_t\co M\to M$ such that $f_t\circ \psi_0=\psi_t.$
\end{lem}

\begin{figure}[ht]
	{\epsfxsize=5in\centerline{\epsfbox{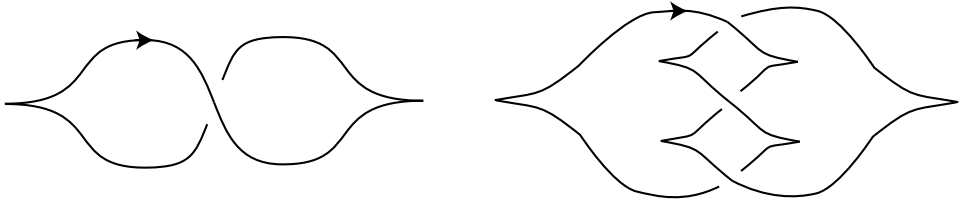}}}
	\caption{Examples of Legendrian knots in the $xy$-projection.}
	\label{fig:examples}
\end{figure}

The contact planes $\xi$ provide a canonical framing of the normal bundle of $\gamma$. If $\gamma$ is
null-homologous, then we may use the Seifert surface $\Sigma$ to assign an integer to the
canonical framing which we call the \dfn{Thurston-Bennequin invariant} of $\gamma$ and
denote by $\tb(\gamma).$ If we orient $\gamma$, then the oriented unit tangent vector
to $\gamma$ will provide a section of $\xi\vert_\gamma.$ The Euler class 
of $\xi\vert_\Sigma$ relative to this section is called the \dfn{rotation number}
of $\gamma$ and is denoted $r(\gamma).$ 
If $\gamma$ is a Legendrian knot in $\mathbb{R}^3$ with the standard contact structure
$\xi=\hbox{ker}(dz+xdy)$, then
its projection onto the $yz$-plane, called the {\it front projection},
will have no vertical tangencies and at each crossing the arc with the
smallest slope will lie over the other arc. See Figure~\ref{fig:examples}.
Moreover any knot diagram satisfying these conditions will represent a Legendrian knot. In this projection
it is easy to compute the Thurston-Bennequin invariant and the rotation number of $\gamma.$ They are given
by
\begin{equation}\label{tbformula}
	\tb(\gamma)= w(\gamma) - r_c,
\end{equation}
and
\begin{equation}\label{rformula}
	r(\gamma)= \frac{1}{2}(D_c -U_c),
\end{equation}
where $w(\gamma)$ is the writhe of $\gamma,$ $r_c$ is the number of right cusps in the projection,
$D_c$ is the number of downward cusps and $U_c$ is the number of upward cusps in the projection. In 
Figure~\ref{fig:examples}, the unknot has $\tb=-2$ and $r=-1$ and the trefoil has $\tb=-6$ and $r=1$.
It is an interesting exercise to work these formulas out or see \cite{Gompf98}.

\begin{figure}[ht]
	{\epsfxsize=4in\centerline{\epsfbox{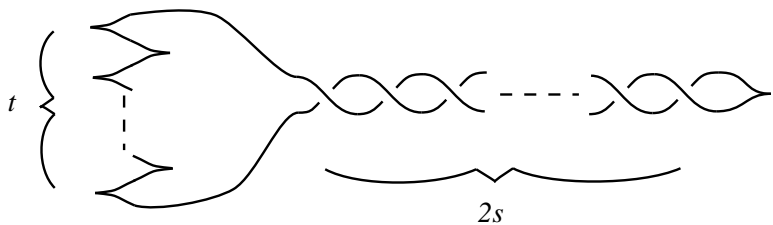}}}
	\caption{Legendrian unknots with $\tb=-2s-t$ and $r=\pm (t-1)$.}
	\label{fig:unknots}
\end{figure}

The only result previously known concerning
the classification of Legendrian knot is due to Eliashberg and Fraser \cite{ef}:
\begin{thm}  \label{unknot}
	Two oriented Legendrian unknots in a tight contact structure 
	are Legendrian isotopic if and only if the
	have the same Thurston-Bennequin invariant and the same rotation
	number. Moreover, a complete list of Legendrian unknots is given in Figure~\ref{fig:unknots}.
\end{thm}
We will present a short proof of this theorem in Section \ref{unknotproof}.

Call a knot type $\mathcal{K}$ \dfn{Legendrian simple} if (oriented) Legendrian knots in this knot type are determined by
their Thurston-Bennequin invariant and rotation number.  Thus the unknot is Legendrian simple.

\begin{figure}[ht]
	{\epsfxsize=3.7in\centerline{\epsfbox{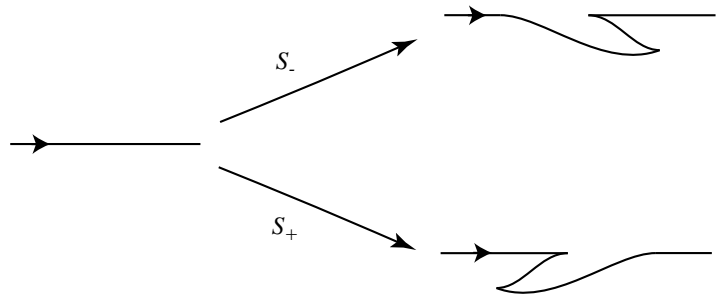}}}
	\caption{Stabilization of Legendrian knots (in $yz$-projection).}
	\label{fig:stabilize}
\end{figure}

\begin{figure}[ht]
	{\epsfxsize=2.5in\centerline{\epsfbox{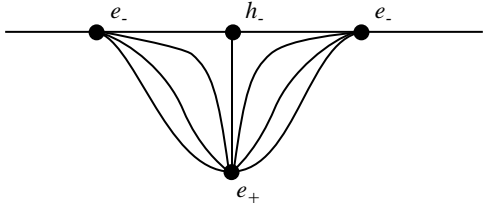}}}
	\caption{A bypass.}
	\label{fig:stabilizingdisk}
\end{figure}

The Bennequin inequality
\begin{equation}\label{LegendrianBennequin}
	\tb(\gamma)+\vert r(\gamma)\vert\leq -\chi(\Sigma)
\end{equation}
provides a (non-optimal, see Theorem~\ref{maxtb}) upper bound on the Thurston-Bennequin invariant
of a Legendrian knot. It is easy, however, to decrease the Thurston-Bennequin invariant. Specifically, 
let $\gamma$ be an oriented Legendrian knot. We can find a contactomorphism from a neighborhood $N$
of $\gamma$ to $M=\{(x,y,z)\in \mathbb{R}^3_y\vert x^2+z^2<\epsilon\},$ where $\mathbb{R}^3_y$ is $\mathbb{R}^3$ 
modulo $y\mapsto y+1$ and $\xi_0=\{dz+xdy=0\},$ 
and $\gamma$ is sent to the image of the $y$-axis in $M$. Now a \dfn{positive (negative)
stabilization} of $\gamma,$ $S_+(\gamma)$ ($S_-(\gamma)$) is the curve in $N$ corresponding to the curve in 
$M$ shown in Figure~\ref{fig:stabilize}. 
From Equations~\ref{tbformula} and \ref{rformula} it is easy to see that 
\begin{equation}
	\tb(S_\pm(\gamma)) = \tb(\gamma)-1
\end{equation}
and
\begin{equation}
	r(S_\pm(\gamma)) = r(\gamma)\pm 1.
\end{equation}
It is important to notice that $S_+(\gamma)\cup \gamma$ cobound a disk $D$ for which $\gamma\setminus\partial D
=S_+(\gamma)\setminus \partial D,$ $\gamma\cap D$ contains three negative singularities, two elliptic and
one hyperbolic, and $S_+(\gamma)\cap D$ contains the same two elliptic singularities and a positive
elliptic singularity. See Figure~\ref{fig:stabilizingdisk}.
We also observe the following simple lemma (cf. \cite{ft}).
\begin{lem}\label{aboutstabilization}
	Stabilization is well-defined and $S_+(S_-(K))=S_-(S_+(K)).$
\end{lem}

In this paper we use a strategy developed in \cite{Etnyre99} (for transversal knots) to classify various 
Legendrian knots. Specifically:
\begin{strategy}\label{mainstrat}
	Given a knot type $\mathcal{K}$,
	\begin{enumerate}
		\item Classify Legendrian knots realizing $\mathcal{K}$ with maximal
			Thurston-Ben\-ne\-quin invariant.
		\item Show that all Legendrian knots realizing $\mathcal{K}$ without maximal Thurston-Bennequin
			invariant destabilize ({\em i.e.}, are stabilizations of other Legendrian knots).
		\item If $K$ and $K'$ are Legendrian knots realizing $\mathcal{K}$ with maximal 
			Thurston-Bennequin invariant, then understand the relationship between their stabilizations.
	\end{enumerate}
\end{strategy}
From Lemma~\ref{aboutstabilization}, if all three steps can be carried out, then we will have classified
Legendrian knots realizing $\mathcal{K}.$ 
Note that if $\mathcal{K}$ has a unique Legendrian realization with maximal Thurston-Bennequin invariant 
(as is the case with unknots, positive torus knots, and figure eight knots), then Step 3 is unnecessary.

\subsection{Transversal Knots}

A knot $\gamma$ in a contact manifold $(M,\xi)$ is called \dfn{transversal} if it is everywhere transverse
to the contact planes. One would like to classify transversal knots up to 
isotopy through transversal knots. The analog to Lemma~\ref{extend} for transversal knots says that this is 
equivalent to the classification of transversal knots up to global contact isotopy.

\begin{figure}[ht]
	{\epsfxsize=3.5in\centerline{\epsfbox{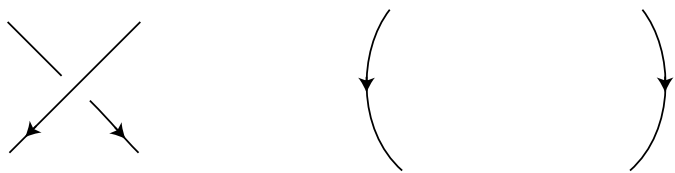}}}
	\caption{Diagram segments not allowed in $yz$-projection of a transversal knot.}
	\label{fig:transversal}
\end{figure}

In addition to the topological knot type, (nullhomologous) transversal knots have one other classical invariant --- the
self-linking number. If $\Sigma$ is the Seifert surface for the transversal knot $\gamma,$ then we can find a nonzero
section $\sigma$ of $\xi\vert_\Sigma.$ We may use $\sigma$ to push off a parallel copy $\gamma'$ of $\gamma$ and then
define the \dfn{self-linking number} of $\gamma$ to be 
\begin{equation}
	l(\gamma)= \gamma\cdot\Sigma,
\end{equation}
where $\cdot$ denotes oriented intersection number ({\em i.e.}, $l(\gamma)$ is just the linking number of $\gamma$ and
$\gamma'$). 
Note if $\xi$ and $M$ are oriented, then $\gamma$ has a natural induced orientation.
If $\gamma$ is a (generic) transversal knot in $\mathbb{R}^3$ with the standard contact structure
$\xi=\hbox{ker}(dz+xdy)$, then
its projection onto the $yz$-plane is a diagram in which one does not see the diagram segments in
Figure~\ref{fig:transversal}.
Moreover, any knot diagram not containing these diagram segments will represent a transversal knot. In this projection
it is easy to compute the self-linking number of $\gamma.$ It is given by
\begin{equation}
	l(\gamma)=w(\gamma),
\end{equation}
where $w(\gamma)$ is the writhe of $\gamma.$

Until recently, the only classification results known for transversal knots were the following two results:
\begin{thm}[Eliashberg \cite{el:knots}]
	Two transversal unknots are trans\-ver\-sal\-ly isotopic if and only if their self-linking numbers agree.
	Moreover, the self-linking numbers of transversal unknots are precisely the negative odd integers.
\end{thm}

\begin{thm}[Etnyre \cite{Etnyre99}]
	Two trans\-ver\-sal positive torus \-knots are trans\-ver\-sal\-ly isotopic if and only if their self-linking
	numbers agree and they are of the same topological knot type. Moreover, the self-linking numbers of
	a transversal positive $(p,q)$-torus knot ($p,q>0$) are precisely the odd integers less than or equal to 
	$pq-p-q.$
\end{thm}
Following Birman and Wrinkle \cite{bw}, one says a knot type $\mathcal{K}$ is \dfn{transversally simple} if
transversal knots in this knot type are determined by their self-linking number.
Recently, using an interesting connection to braid theory, the above theorem was greatly generalized:
\begin{thm}[Birman and Wrinkle \cite{bw}, Menasco \cite{M}]
	Iterated torus knots are trans\-ver\-sally simple.
\end{thm}
We recall the definition of iterated torus knots. If $\gamma$ is any knot, let $N=S^1\times D^2$ be a neighborhood of
$\gamma$ so that $S^1\times \{p\},$ where $p\in\partial D^2,$  bounds a Seifert surface in the complement of $N.$
If the integers $p$ and $q$ are relatively prime, then the homology class $p(\partial D^2) + q (S^1\times\{p\})$ can
be represented by an embedded closed curve $\gamma_{(p,q)}\subset \partial N.$ The curve $\gamma_{(p,q)}$ is called
the \dfn{$(p,q)$-cable} of $\gamma.$ Any iterated cable of the unknot is called an \dfn{iterated torus knot}.
The above theorem reduces the transversal classification of
iterated torus knots to an existence question, i.e., which self-linking numbers are actually realized by transversal
knots in these knot types. In private communication with the authors, Birman has finished the classification by showing that the
possible self-linking numbers for $(\ldots(\gamma_{(p_1,q_1)})_{(p_2,q_2)}\ldots)_{(p_n,q_n)},$ where $\gamma$ is 
an unknot, are odd integers 
less than or equal to $l_n$.  Here $l_n$ is computed as follows: Let $b_n=q_1q_2\ldots q_n,$ where we use
the convention that $0<q_i<\vert p_i\vert,$  and set $a_1 = (q_1-1)p_1$ and $a_i=(q_i-1)p_i -a_{i-1}q_i^2.$
Then $l_n=a_n-b_n.$ 

The Bennequin inequality for transversal knots is
\begin{equation}\label{TransversalBennequin}
	l(\gamma)\leq -\chi(\Sigma),
\end{equation}
where $\Sigma$ is a Seifert surface for $\gamma.$ The Bennequin inequality for Legendrian knots, 
Equation~\ref{LegendrianBennequin}, can be obtained from this inequality as follows: If $\gamma'$ is a Legendrian
knot, then we may embed an annulus $A=\gamma'\times[-\epsilon, \epsilon]$ in a neighborhood of
$\gamma'$ so that $\gamma'$ is the core, $\gamma'\times\{0\},$ of $A$ and
the framing $A$ induces on $\gamma'$ is the same as the one induces by $\xi.$ Now if $A$ is thin enough we can 
assume the characteristic foliation contains $\gamma'$ as a closed leaf and the other leaves are transverse to
the boundary and spiral to $\gamma'.$  Now $\gamma'\times\{\epsilon\}$ is a positive
transversal unknot which we denote $T_+(\gamma')$ and call the \dfn{positive transversal push-off} of $\gamma'.$
Similarly we have the \dfn{negative transversal push-off,} $T_-(\gamma'),$ of $\gamma'.$ One can show (\cite{el:twenty})
\begin{equation}
	l(T_\pm(\gamma'))=\tb(\gamma')\pm r(\gamma'),
\end{equation} 
and Equation~\ref{LegendrianBennequin} follows from this and Equation~\ref{TransversalBennequin}.

One can define a notion of stabilization for transversal knots, as we did for Legendrian one, and use this
to develop a strategy to classify transversal knots. In fact this was used in \cite{Etnyre99} to classify transversal
positive torus knots. In this paper we use a different approach originally described to the authors
by Eliashberg. We begin with a definition. Two Legendrian knots $\gamma$ and $\gamma'$ are called \dfn{stably isotopic} if
there is some $n$ and $n'$ such that $S^n_+(\gamma)$ and $S^{n'}_+(\gamma')$ are Legendrian isotopic. 
Since $\tb(\gamma)+r(\gamma)=\tb(S_+(\gamma))+r(S_+(\gamma)),$ a natural invariant of the stable isotopy 
class of a Legendrian knot is 
\begin{equation}
	s(\gamma)=\tb(\gamma)+r(\gamma),
\end{equation}
which we call the \dfn{stable Bennequin invariant.}
We call a knot type $\mathcal{K}$ \dfn{stably simple} if Legendrian knots in this knot type are stably isotopic if and
only if their stable Bennequin invariants agree.
Note that we allow ourselves to use only {\it positive} stabilizations in the definition of
{\it stable isotopy}. The reason for this is twofold: first, if we allow negative
stabilizations as well, then any two topologically isotopic Legendrian knots would be stably isotopic (cf. \cite{ft}).
And, second, we have the
following:  
\begin{thm}\label{thm:stably-transversally}
	A knot type $\mathcal{K}$ is stably simple if and only if it is transversally simple.
\end{thm}

\begin{proof}
Begin by assuming $\mathcal{K}$ is stably simple.
Let $\gamma$ and $\gamma'$ be two transversal knots in the knot type $\mathcal{K}$ with the
same self-linking numbers. Note $\gamma$ has a neighborhood
$N$ contactomorphic to $\{(r,\theta, z)\in \mathbb{R}^2\times S^1\vert r<\epsilon\}$ with the contact structure 
given by $\{(\alpha=dz+r^2d\theta)=0\}.$ Now, for large integers $n,$ if $T_n$ are the tori in $N$ with
$r=\frac{1}{\sqrt{n}}$, then the characteristic
foliation on $T_n$ is by $(-1,n)$ curves. Let $L_n$ be a leaf in this characteristic foliation. Note $L_n$ is a 
Legendrian knot topologically
isotopic to $\gamma$ and that $S_+(L_n)=L_{n+1}$. Note also that if we have some $L_n$, then we have
$L_m$ for all $m\geq n.$ Thus from $\gamma$ we have a well-defined stable isotopy class of Legendrian knot $L_n$ and
from $\gamma'$ we similarly get $L'_n.$ Now since $s(L'_n)=l(\gamma')=l(\gamma)=s(L_n)$ we know that there is 
some $m>n$ and $m'>n'$ such that $L_m$ is Legendrian isotopic to $L'_{m'}.$ Finally, by observing that $T_+(L_m)=\gamma$
and $T_+(L'_{m'})=\gamma'$, we see that $\gamma$ and $\gamma'$ are transversally isotopic. Thus $\mathcal{K}$ is
transversally simple.

The other implication is proved at the end of the next section as it requires the theory of convex surfaces developed 
there.
\end{proof}

\section{Convexity in contact geometry}

\subsection{Contact vector fields}

A vector field $v$ is called a \dfn{contact vector field} on a contact manifold $(M,\xi)$
if the flow of $v$ preserves the contact structure $\xi$. 
A surface
$\Sigma$ in $M$ is called \dfn{convex} if there is a contact vector field
transverse to $\Sigma$.  Generically, $v\vert_\Sigma$ will be tangent to $\xi\vert_\Sigma$ along curves $\Gamma$
that divide $\Sigma$ in a special way. In general, consider
$\mathcal{F}$ a singular foliation on an orientable surface $\Sigma$ and $\Gamma$ a disjoint union 
of simple closed curves on $\Sigma.$ We say $\Gamma$ \dfn{divides} $\mathcal{F}$ if $\Gamma$ is
transverse to $\mathcal{F},$ $\Sigma\setminus\Gamma$ is the disjoint union of two 
(possibly disconnected) surfaces $\Sigma_+$ and $\Sigma_-$ 
with $\partial\overline{\Sigma_+}=-\partial\overline{\Sigma_-}=\Gamma,$ 
and there is a vector field $u$ and
volume form $\omega$ on $\Sigma$
so that $u$ is tangent to $\mathcal{F},$ $\pm \mathcal{L}_u\omega>0$ on $\Sigma_\pm$ and $u\vert_\Gamma$
points out of $\Sigma_+.$ 
We refer the reader to \cite{gi:convex} and \cite{k} for
proofs of the following facts: 
\begin{itemize}
\item Contact vector fields may be identified with a section of the bundle $TM/\xi.$ 
\item A closed surface may be isotoped by a $C^\infty$-small isotopy so that it is convex.
\item Let $\Sigma$ be a compact orientable surface in the contact manifold $(M^3,\xi)$.  Then 
	$\Sigma$ is a convex surface if and  only if there is a tubular neighborhood $N$ of
	$\Sigma$ in $M$ that is contactomorphic to  
	$(\Sigma\times(-\epsilon,\epsilon),\beta+u\,dt)$ taking
	$\Sigma$ to $\Sigma\times\{0\}$, where
	$\beta$ is a 1-form on $\Sigma$ and $u$ is a function on $\Sigma$ and $\epsilon>0$.
\item Let $\Sigma$ be a closed orientable surface. Then $\Sigma$ is convex if and only if
	there is a dividing set $\Gamma$ for $\Sigma_\xi.$ (If there is any ambiguity, we will also write 
	$\Gamma_\Sigma$ instead of $\Gamma$ to denote a dividing set of $\Sigma$.)
	Two dividing sets $\Gamma$ and $\Gamma'$ for the same $\mathcal{F}$ are isotopic; hence we
	will slightly abuse language and refer to $\Gamma$ as `the' dividing set of $\Sigma$.
\end{itemize}

Let $v$ be a contact vector field for $(M,\xi)$ that is transverse to a surface $\Sigma$ and
let $\Gamma$ be the dividing curves on $\Sigma.$ 
An isotopy $F:\Sigma\times[0,1]\to M$ of $\Sigma$ is called \dfn{admissible} if $F(\Sigma\times\{t\})$
is transversal to $v$ for all $t.$ The major result concerning convex surfaces says that up to admissible
isotopies the dividing set dictates the geometry of $\xi$ near $\Sigma.$ More specifically, we have the following
Flexibility Theorem:
\begin{thm}[Giroux \cite{gi:convex}, Kanda \cite{K2}]\label{giroux:main}
	Let $\Sigma$ be a closed surface or a surface with Legendrian boundary.
	Let $\Gamma$ be the dividing set for $\Sigma_\xi$ and $\mathcal{F}$ another
	singular foliation on $\Sigma$ divided by $\Gamma.$ Then there is an admissible 
	isotopy $F$ of $\Sigma$ such that $F(\Sigma\times\{0\})=\Sigma,$ 
	$F(\Sigma\times\{1\})_\xi=\mathcal{F}$ and the isotopy is fixed on $\Gamma.$
\end{thm}

\begin{rem}
\em{It is useful to note, and we will frequently implicitly use, that in a tight structure no dividing curve
on a convex closed surface or a convex surface with Legendrian boundary can bound a disk unless the surface
is a 2-sphere. This observation is due to Giroux.  For example, see \cite{Honda1}.}
\end{rem}

A useful formulation of Giroux's Flexibility Theorem above is called the Legendrian 
Realization Principle, which is due to Kanda \cite{K2}:

\begin{thm}[Legendrian Realization Principle] \label{lerp}
Consider $C$, a closed curve, on a closed convex surface or a convex surface $\Sigma$
with Legendrian boundary.  Assume $C\pitchfork \Gamma_\Sigma$ and every component of
$\Sigma\backslash C$ nontrivially intersects $\Sigma$.
Then there exists an
admissible isotopy $F_t=F(\cdot, t)$, $t\in[0,1]$ so that
\be
\item $F_0=id$,
\item $F_t(\Sigma)$ are all convex,
\item $F_1(\Gamma_\Sigma)=\Gamma_{F_1(\Sigma)}$,
\item $F_1(C)$ is Legendrian.
\ee
\end{thm}

If $\gamma$ is a closed oriented Legendrian curve in a surface $\Sigma$, then we define the \dfn{twist of $\gamma$ relative
to $\Sigma$}, $t_\Sigma(\gamma),$ to be the twisting of the contact planes $\xi$ along $\gamma$ measured with respect 
to the framing induced by $\Sigma.$ Note that if $\Sigma$ is a Seifert surface for $\gamma$, then
$t_\Sigma(\gamma)=\tb(\gamma).$ 
We now have the following useful theorem (cf. \cite{K2} as well as \cite{ef} for the first sentence of the theorem):
\begin{thm}[Kanda \cite{K2}]\label{makeconvexrelcurve}
	If $\gamma$ is a Legendrian curve in a surface $\Sigma,$ then $\Sigma$ may be isotoped relative to
	$\gamma$ so that it is convex
	if and only if $t_\Sigma(\gamma)\leq 0.$ Moreover, if $\Sigma$ is convex, then
	\begin{equation}
		t_\Sigma(\gamma)= -\frac{1}{2}\#(\gamma\cap\Gamma),
	\end{equation}
	where $\Gamma$ is the set of dividing curves for $\Sigma_\xi.$
\end{thm}

Here, $\#(a\cap b)$ is the cardinality of $a\cap b$.  Note that $\gamma\pitchfork\Gamma$ because
$\Gamma\pitchfork \Sigma_\xi$.

\subsection{Surfaces with boundary and convexity}

Let $(M,\xi)$ be a contact manifold.
Throughout this section let $\Sigma$ be a compact surface whose boundary $\bdry\Sigma$ is
everywhere tangent
to $\xi.$ 
Moreover, whenever we are considering oriented Legendrian curves we assume $\Sigma$ is oriented and the orientation is
consistent with the orientation of $\partial\Sigma.$ 
From Theorem~\ref{makeconvexrelcurve} we have 
\begin{lem}
	The surface $\Sigma$ may be made convex if and only if the twist of
	$\xi$ about each boundary component is less than or equal to zero.
\end{lem}
For Legendrian knots we have
\begin{lem}[Kanda \cite{K2}]\label{lem:invts}
	Suppose $\Sigma$ has a single boundary component $\gamma,$ and $\gamma$ is Legendrian.
	Then $\Sigma$ may be made convex
	if and only if $\tb(\gamma)\leq 0.$ Moreover, if $\Sigma$ is convex with dividing curves
	$\Gamma$, then
	\begin{equation}
		\tb(\gamma)=-\frac{1}{2}\# (\gamma\cap\Gamma)
	\end{equation}
	and
	\begin{equation}\label{formulaforrfromconvexity}
		r(\gamma)=\chi(\Sigma_+)-\chi(\Sigma_-),
	\end{equation}
	where $\Sigma_\pm$ are as in the definition of convexity.
\end{lem}

We will also need the following Edge-Rounding Lemma where two convex surfaces meet `perpendicularly'
along a common Legendrian boundary. But first we discuss a normal form for two intersecting
convex surfaces.
Let $\Sigma_1$ and $\Sigma_2$
be compact convex surfaces with Legendrian boundary in a contact manifold
$(M,\xi)$, $\Sigma_1$ and $\Sigma_2$ intersect transversally along a common boundary
component $\gamma$, and $t_{\Sigma_1}(\gamma)=t_{\Sigma_2}(\gamma)<0$.  Then
$\Sigma_1$, $\Sigma_2$ can be $C^0$-small perturbed, fixing the boundary, so that
the neighborhood $N(\gamma)$ of $\gamma$ is locally contactomorphic to 
$\{(x,y,z)\in\R_z^3\vert x^2+y^2<\epsilon\}$ with the contact structure $\alpha_n=\sin(2\pi nz)dx+\cos(2\pi nz)dy$
for some $n\in \Z^+,$ 
where $\R_z^3$ is $\R^3$ modulo $z\mapsto z+1.$
$\Sigma_1\cap N(\gamma) = \{x = 0, 0 \le
y \le \varepsilon\}$
and $\Sigma_2\cap N(\gamma) = \{y = 0, 0 \le x \le \epsilon\}$. 
\bbr\label{transfer}
This local model for the intersection
of two convex surfaces shows how to transfer information from one convex surface to another.
\eer

\begin{lem}[Edge-Rounding] \label{edge}   In the situation described above, if we join $\Sigma_1$
and $\Sigma_2$ along
$x = y = 0$ and round the common edge, the resulting surface is convex, and the
dividing curve
$z = {k\over 2n}$ on $\Sigma_1$ will connect to the dividing curve $z = {k
\over 2n} - {1\over
4n}$ on $\Sigma_2$, where $k = 0, \cdots, 2n - 1$. Here we are assuming that $\Sigma_1$ and
$\Sigma_2$ have been oriented so that the orientations agree after rounding.
\end{lem}

\subsection{Bypasses and simplifications of the characteristic foliation}

We now wish to describe a method for altering the dividing curves of a convex surface.
These techniques were first exploited in \cite{Honda1} and precursors to them appeared
in \cite{k}.

\begin{figure}[ht]
	{\epsfxsize=2.5in\centerline{\epsfbox{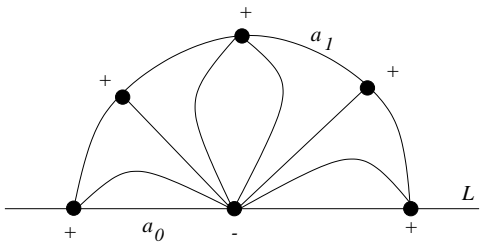}}}
	\caption{A bypass.}
	\label{fig:bypassdisk}
\end{figure}

Given a Legendrian curve $\gamma$ a \dfn{bypass for $\gamma$} (see Figure~\ref{fig:bypassdisk}) 
is a convex disk $D$ with
boundary made up of two smooth Legendrian arcs $a_0$ and $a_1$ where: 
\begin{itemize}
	\item $a_0=(\gamma\cap D)\subset \gamma,$ 
	\item along $a_0$ there are three elliptic singularities in $D_\xi$ --- two with 
		the same sign occurring at the endpoints and one with a different sign in the 
		interior of $a_0,$
	\item along $a_1$ all the singularities have the same sign, their local degrees alternate,
		 and there are at least three of them,
	\item there are no interior singularities in $D_\xi.$
\end{itemize}
Note all the singularities of $D_\xi$ have the same sign except one --- the sign of this
singularity is \dfn{the sign of the bypass}.  Note the relation between the definition of
a bypass and a stabilizing disk. In fact, suppose $D'$ is a stabilizing disk for a Legendrian knot
$\gamma'$ and, once stabilized, we obtain the knot $\gamma$.  Then from $\gamma$'s perspective $D'$
is a bypass --- {\em i.e.,\ } $D'$ shows how to isotop $\gamma$ so that it twists less. Moreover if
$D$ is a bypass for a knot $\gamma$ and $\gamma'$ is the knot obtained from pushing $\gamma$ across
$D$, then $D$ may be isotoped (canceling extra singularities along $\alpha_1$) so that it is a
stabilizing disk for $\gamma'.$

Dividing curves are quite helpful in locating bypasses. To see this let $\Sigma$ be a convex surface with
Legendrian boundary. In this situation one may manipulate the the characteristic foliation 
to make all the singularities along the boundary (half)-elliptic (\cite{ef}). Now if $\tb(\partial\Sigma)=-n\leq0$
then the dividing curves intersect $\partial\Sigma,$ $2n$ times. Suppose one of these dividing curves is
boundary-parallel --- {\em i.e.,} cuts off a half-disk which has no other intersections with $\Gamma_\Sigma$.
If we think of the characteristic foliation as given by a flow, then we may flow this
dividing curve ``away'' from the boundary and it will limit to a Legendrian curve $\alpha.$ The curve $\alpha$
will separate a disk $D$ from $\Sigma$ which will be a bypass for $\partial\Sigma.$   (If $t(\bdry \Sigma)=-1$,
$\Sigma\not=D^2$, and $\Gamma_\Sigma$ is a single boundary-parallel curve along $\bdry \Sigma$, then
we need to exercise a little care.  In this case, we first need to perturb $\Gamma$ and modify $\Gamma_\Sigma$,
and then use a more general form of the Legendrian Realization Principle.  We remark that this case actually
does not appear in this paper.)

\begin{lem}  \label{createbypass}
	Given a boundary-parallel dividing curve $\delta$ on a convex surface $\Sigma$
	with Legendrian boundary, one may find a
	bypass for the boundary, provided $\Gamma_\Sigma$ is not a single arc on $\Sigma=D^2$.
\end{lem}

\begin{rem}
	\em{We will frequently abuse terminology and refer to a boundary-parallel dividing curve as a bypass.}
\end{rem}

If $\Sigma$ is an annulus, we may find boundary-parallel dividing curves using the following
Imbalance Principle (see \cite{Honda1}):

\begin{prop}[The Imbalance Principle]   \label{imba} If $\Sigma = S^1
\times [0,
1]$ is convex and has Legendrian boundary where $t(S^1 \times \{0\}) < t(S^1 \times
\{1\}) \le 0$,
then there exists a boundary-parallel dividing curve (and hence a bypass) along $S^1
\times \{0\}$.
\end{prop}

\begin{figure}[ht]
	{\epsfxsize=4in\centerline{\epsfbox{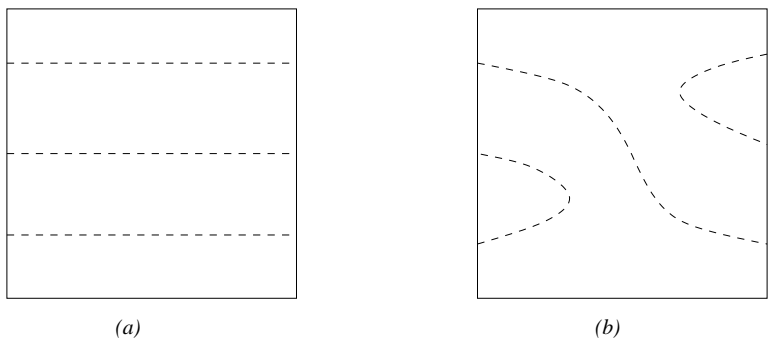}}}
	\caption{Dividing curves on $A$ before the isotopy (a) and after (b).}
	\label{fig:bypassmove}
\end{figure}

Bypasses may be used to alter the dividing curves of a convex surface as follows:
\begin{prop}\label{prop:bypass}
	Let $A=[0,1]\times [0,1]$ be a convex square with three horizontal dividing curves and
	vertical ruling. Let $\gamma$ be one of the vertical ruling curves and $D$ a bypass
	for $\gamma$ disjoint from $A.$ Then we may isotop $A$ rel boundary by pushing $A$ across $D$
	so as to alter the characteristic foliation as shown in Figure~\ref{fig:bypassmove}. 
\end{prop}

\subsection{The classification of tight contact structures on solid tori}\label{torusstuff}

Let $T$ be  convex torus. Assuming $\xi$ is tight, we know that no dividing curve bounds a disk, and
hence the dividing curves are parallel essential curves.
If there are $2n$ parallel curves, then, using Theorem~\ref{giroux:main}, we
may assume there are $2n$ curves of singularities in $T_\xi,$ one in each region of the complement of the dividing
curves. We call these curves the \dfn{Legendrian divides}, and their slope will be called the \dfn{boundary slope}.
We may also assume the the other leaves in $T_\xi$ form a 1-parameter family of closed curves. We call these
curves the \dfn{Legendrian ruling curves}. It is important to notice that we can make the slope of these
ruling curves whatever we wish, except the boundary slope.
If the characteristic foliation on a  convex torus has this nongeneric 
form we say the $T_\xi$ is in \dfn{standard from}. If $n>1$ and we can find a bypass for one of the ruling
curves, then we can isotop $T$ so as to reduce $n$ by one. If $n=1$ we have:
\begin{thm}[Honda \cite{Honda1}]\label{bypassalteration}
	Let $T$ be a convex torus in standard from with two dividing curves and in some basis for $T$ the slope of the 
	dividing curves is 0. If we find a bypass on a ruling curve of slope between $-\frac{1}{m}$ and $-\frac{1}{m+1},$
	$m\in\Z,$ then after pushing $T$ across the bypass the new torus has two dividing curves of slope $-\frac{1}{m+1}.$
\end{thm}
This theorem was a key step in proving:
\begin{thm}[Honda \cite{Honda1} (see also \cite{Giroux99})]\label{toriclassification}
	Let $(r_0,\ldots,r_k)$ give the continued fraction decomposition of $-\frac{p}{q},$ where $p>q>0.$
	Then there are $\vert (r_0+1)\ldots(r_{k-q}+1)(r_k)\vert$ tight contact structures on $S^1\times D^2$
	with standard convex boundary having two Legendrian divides of slope $-\frac{p}{q}.$ (Here our convention
	is that the meridian has slope $0$.) Moreover, all these
	structures are distinguished by the number of positive regions on a convex meridional disk
	with Legendrian boundary.
\end{thm}

It is not hard to see that if we have the `standard' universally tight contact structure on $T^2\times I$  (with
coordinates $(x,y,z)$)
given by $\alpha=\sin({\pi\over 2} z)dx+\cos({\pi\over 2} z)dy$, and we take the
one-parameter family of tori $T^2\times \{z\}$ with linear characteristic foliations, then
the slopes of these foliations must decrease as $z$ increases.
For convex tori this is a little more complicated but still true.
\begin{lem}\label{torustwisting}
	If $T^2\times [0,1]$ has convex boundary in standard form and the boundary slope on $T^2\times \{i\}$ 
	is  $s_i$, for $i=0,1,$ then we can find 
	convex tori parallel to $T^2\times \{i\}$ with any boundary slope $s$  in $[s_1, s_0]$ 
	(if $s_0<s_1$ then this means $[s_1,\infty]\cup [-\infty, s_0]$).   
\end{lem}
This follows from the classification of tight
contact structures on $T^2\times [0,1]$  (see \cite{Honda1, Giroux99, EH}).
From this Lemma one can easily show:
\begin{lem}\label{solidtorustwisting}
	If $S=D^2\times S^1$ has convex boundary with boundary slope $s<0,$ then we can find a convex torus
	parallel to the boundary of $S$ with any boundary slope in $[s,0).$
\end{lem}

We also will make use of the following consequence of the classification of tight
contact structures on $T^2\times I$:
\begin{lem}\label{model} Consider a tight contact structure on $T^2\times [0,1]$ with boundary slopes
$s_1=-{1\over m}$, $s_0=-{1\over m+1}$ ($m\in \Z^+$) (for $T^2\times\{1\}$ and $T^2\times \{0\}$, respectively).
If $s_1<s<s_0$, then there exists a pre-Lagrangian (= linearly foliated) torus $T$ parallel to $T^2\times\{i\}$, and
every convex surface $T'$ in standard form with slope $s$, after contact isotopy, is transverse to $T$, and
$T\cap T'$ is exactly the union of the Legendrian divides of $T'$.
\end{lem}

We shall also have need for the relative Euler class of a contact structure on $U=T^2\times [0,1].$ Following
\cite{EH}, if $\xi$ is a tight contact structure on $U$ for which the boundary is convex and in standard form, then
let $v$ be a nonzero vector field transverse to and twisting (with $\xi$) along the Legendrian ruling curves
and tangent to the Legendrian divides. Now let $e$ be the relative Euler class in $H^2(U,\partial U;\Z)$
for this section of the bundle $\xi\vert_{\partial U}.$ It is important to note that  $e$ is unchanged if we perform a 
$C^0$-small isotopy of $\partial U$ so as to alter the slopes of the ruling curves. Now if $c$ is an
oriented curve on $T^2\times\{0\}$, then we can assume the annulus $A=c\times [0,1]$ has Legendrian
boundary and is convex. A slight generalization of Equation~\ref{formulaforrfromconvexity} yields
\begin{equation}\label{formulafore}
	e(c)\equiv e(A)=\chi(A_+)-\chi(A_-),
\end{equation}
where $A_\pm$ are the positive and negative regions of $A$ from the definition of convex surface.
See \cite{EH, K2} for details.

\subsection{Legendrian Unknots}  \label{unknotproof}

In this section we present a brief proof of Theorem~\ref{unknot}, using Strategy~\ref{mainstrat}.
Consider a contact 3-manifold $(M,\xi)$.  The maximal Thurston-Bennequin invariant
for a knot type $\mathcal{K}$ is an invariant of the knot type, and will be denoted
$\overline{\tb}(\mathcal{K})$.
Observe that the tightness of $\xi$ is equivalent to $\overline{\tb}(\mathcal{K})<0$, where
$\mathcal{K}$ is the knot type for the unknot.  Since there exists a Legendrian unknot
$L$ with $\tb(L)=-1$ and $r(L)=0$, $\overline{\tb}(\mathcal{K})=-1$.  

We begin with Step 2 in Strategy~\ref{mainstrat}.
If $\tb(K')<-1$, then take a spanning convex disk $D$.  There must always exist a
bypass by Lemma \ref{createbypass} --
hence every $K'$ with $\tb(K')<-1$ is a stabilization of a Legendrian unknot with $\tb(K)=-1.$ 

We now complete the proof by showing that there is a unique Legendrian unknot with $\tb(K)=-1.$ 
In this part of the proof we assume, for clarity, that $M=S^3.$ The general case is not much more
difficult but obscures the main ideas.
Take a Legendrian unknot $K$ with
$\tb(K)=-1$ and consider a spanning convex disk $D$ with $\bdry D=K$. 
Since $\tb(K)=-1$, there is only one possible dividing set, and we may  uniquely
normalize the characteristic foliation using the Flexibility Theorem (Theorem \ref{giroux:main}).
If there are two Legendrian unknots $K$, $K'$ with
$\tb=-1$, then there is a diffeomorphism $f:S^3\to S^3$ (taking $D$ to $D'$) which is a contactomorphism when restricted
to neighborhoods $N(D)$, $N(D')$ of the spanning disks $D$, $D'.$  
The map $f$ may now be isotoped (relative to $N(D)$) into a contactomorphism on all of $S^3$ using 
Theorem~\ref{thm:uniqueB3}, since the tight contact structures on the 3-balls $\overline{S^3\backslash N(D)}$ and
$\overline{S^3\backslash N(D')}$ induce the same characteristic foliations on their boundaries.
Thus, using Theorem~\ref{thm:htpyequiv}, we find that the $\tb=-1$ Legendrian unknot inside $M=S^3$
is unique up to contact isotopy.
  
For $M$ arbitrary, we need to use the fact (whose proof we leave to the reader) that two convex 3-balls with
the same boundary characteristic foliation are contact isotopic (this essentially follows from 
Theorem \ref{thm:uniqueB3} as well). We may then reduce to the $M=S^3$ case, thus proving
that every Legendrian unknot is determined by $\tb$ and $r$.

\subsection{Transversal simplicity implies stable simplicity} 

In this section we complete the proof of Theorem~\ref{thm:stably-transversally} by showing transversal
simplicity implies stable simplicity. Begin by assuming that $\mathcal{K}$ is a transversally simple
knot type. Let  $\gamma$ and $\gamma'$ be two Legendrian knots in this knot type with the same
stable Bennequin invariants, $s(\gamma)=s(\gamma').$ Let $K=T_+(\gamma)$ and $K'=T_+(\gamma')$ be 
the positive transversal push-offs of $\gamma$ and $\gamma'.$  We can find solid tori neighborhoods
$N$ and $N'$ of $K$ and $K'$ with convex boundaries such that $\gamma$ and $\gamma'$ lie, respectively, on the boundary of 
$N$ and $N'.$ Note by thickening the annulus from the definition of $T_+(\gamma)$ we may assume that $\gamma$
is a Legendrian divide on $\partial N$ that is also a $(-1,n)$ curve.
Now, since $\mathcal{K}$ is transversally simple and $l(K)=s(\gamma)=s(\gamma')=l(K')$, we
can find a global contact isotopy taking $K'$ to $K.$ Thus we may think of $N$ and $N'$ as solid torus
neighborhoods of the same transversal knot $K.$ 

In $N\cap N'$ we can find a neighborhood $N_\epsilon$
of $K$ contactomorphic to  $\{(r,\theta, z)\in \mathbb{R}^2\times S^1\vert r<\epsilon\}$ with the contact structure 
given by $\{(\alpha=dz+r^2d\theta)=0\}.$ As before for large integers $m,$ if $T_m$ are the tori in $N_\epsilon$ with 
$r=\frac{1}{\sqrt{m}}$, then the characteristic
foliation on $T_m$ is by $(-1,m)$ curves. Fix some large $m$ and let $L_m$ be a leaf in the characteristic foliation.
We now show that there is some integer $k$ such that $S_+^k(\gamma)$ is Legendrian isotopic to $L_m.$ The same
proof will work for $\gamma'$, thus proving they are stably equivalent. Isotop $T_m$ to be convex with $L_m$ as a
Legendrian divide and the ruling curves all in the class $(-1,n).$ Let $A=\gamma\times [0,1]$ be an annulus embedded in
$N\setminus N_\epsilon$ with $S^1\times \{0\}$ a ruling curve on $\partial T_m$ and 
$S^1\times\{1\}=\gamma.$  If we now make $A$ convex, then
clearly $t_A(\gamma)=0$ and $t_A(L_m)<0$ (if $t_A(L_m)=0$ then $\gamma$ 
is Legendrian isotopic to $L_m$ by Lemma~\ref{model}). Now, by the Imbalance Principle (Proposition~\ref{imba}),
there exist dividing curves boundary-parallel to $S^1\times\{0\}$ (and of course non-boundary-parallel
to $S^1\times\{1\}$). This demonstrates that $L_m$ is a stabilization of $\gamma$ --- however, there might be
positive and negative stabilizations. To see that there are only positive stabilizations we note that $T_+(L_m)$ is
a transversal curve isotopic to the core curve in $N_\epsilon.$ Thus $\tb(L_m)+r(L_m)=l(K)=\tb(\gamma)+r(\gamma),$
which of course implies that all the stabilizations were positive.

\section{Torus Knots}

Let $T^2$ be a standardly embedded torus in $S^3.$ By this we mean that $T$ provides 
a genus one Heegaard splitting of $S^3,$ so $S^3=V_0\cup_T V_1$ where the $V_i$ are 
solid tori.  Let $\mu$ be the unique curve on $T$ that bounds a disk in $V_0$ and
$\lambda$ the unique curve that bounds a disk in $V_1.$ Orient $\mu$ arbitrarily and 
then orient $\lambda$ so that $\mu, \lambda$ form a positive basis for $H_1(T)$ where 
$T$ is oriented as the boundary of $V_0.$ Up to homotopy any curve on
$T$ can be written as $p\mu + q\lambda,$ we shall denote this curve by $K_{(p,q)}.$   
If $p$ and $q$ are relatively prime
then $K_{(p,q)}$ is called a \dfn{$(p,q)$-torus knot.}  We will always assume $|p|>q>0$ which 
affords no loss of generality since we may switch the roles of $V_0$ and $V_1$ and 
$K_{(-p,-q)}$ and $K_{(p,q)}$ describe the same knot.  One
may easily compute that the Seifert surface of minimal genus for $K_{(p,q)}$
has Euler number $|p|+|q|-|pq|.$  

\subsection{Legendrian Torus Knots}

Bennequin's inequality for a Legendrian torus knot says
\begin{equation}
	\tb(K_{(p,q)})+\vert r(K_{(p,q)}\vert\leq -|p|-|q|+|pq|.
\end{equation}
Let $\overline{\tb}_{(p,q)}$ the maximal Thurston-Bennequin invariant of a 
Legendrian $K_{(p,q)}$ then 
\begin{thm}\label{maxtb}
	If $p,q>0$ then
	\begin{equation}
		\overline{\tb}_{(p,q)}=pq-p-q,
	\end{equation}
	and if $p<0$, $q>0$ then
	\begin{equation}
		\overline{\tb}_{(p,q)}=pq.
	\end{equation}
\end{thm}

\begin{rem}
	{\em Note that Bennequin's inequality is sharp for positive torus knots but is
	not for negative torus knots. As mentioned in the introduction, the upper bounds 
	for the Thurston-Bennequin invariant coming from the work of Fuchs and Tabachnikov \cite{ft,Ta} 
	are also not sharp when $pq<0$ and $q$ is odd.}
\end{rem}

We now state our main theorem.

\begin{thm}\label{main}
	Legendrian torus knots in the standard tight contact structure on $S^3$ 
	are determined up to Legendrian isotopy by their knot type, Thurston-Bennequin invariant
	and rotation number.
\end{thm}

To complete the classification of Legendrian torus knots we can also show:

\begin{thm}\label{rrange}
	Let $K$ be a Legendrian $(p,q)$-torus knot with maximal Thurston-Bennequin invariant.
	If $p,q>0,$ then $r(K)=0.$ If $p<0$, $q>0,$ then
	$$r(K)\in\{ \pm(\vert p\vert-\vert q\vert - 2qk) : k\in\mathbb{Z},
	0\leq k<  \frac{\vert p\vert-\vert q\vert}{\vert q\vert}\}.$$
\end{thm}

\begin{figure}[ht]
	{\epsfxsize=3.75in\centerline{\epsfbox{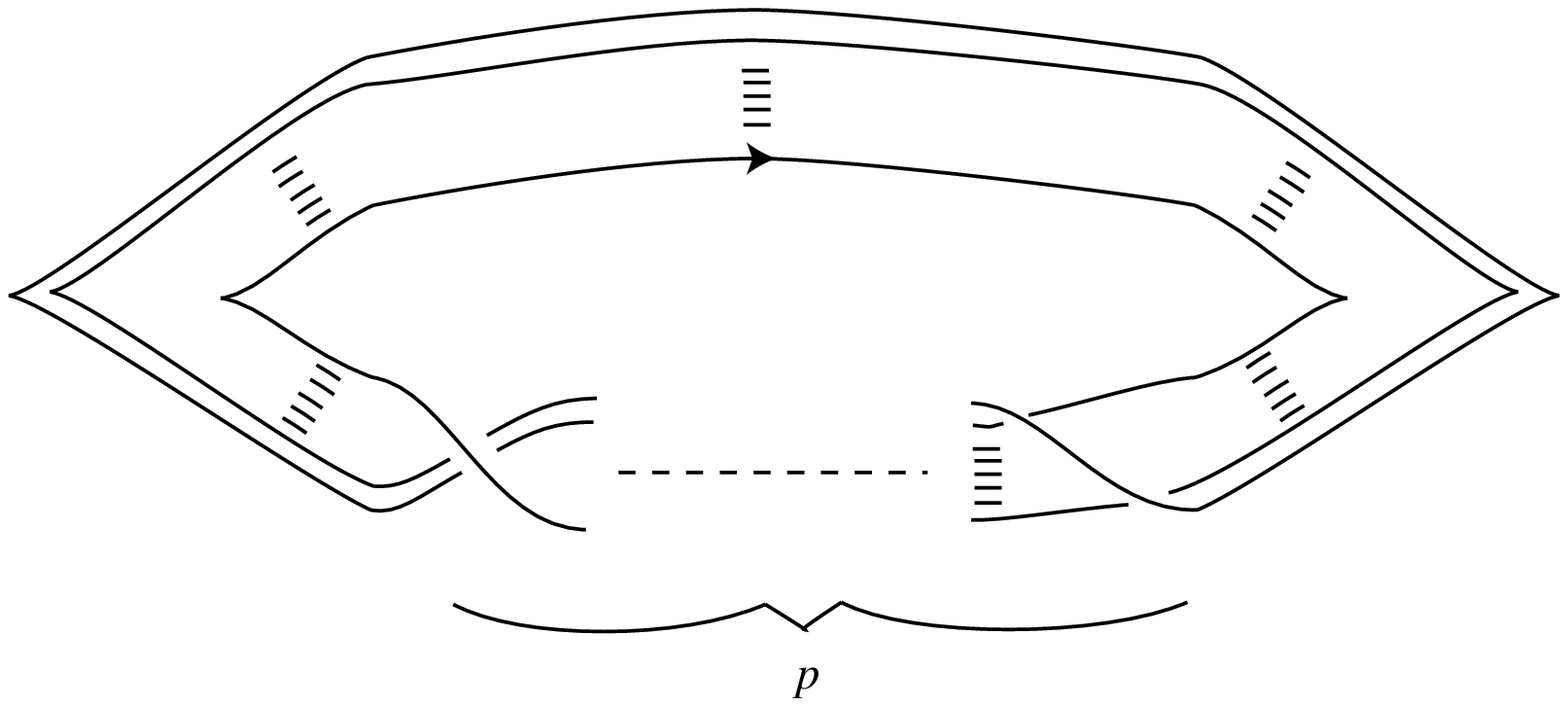}}}
	\bigskip
	\bigskip
	{\epsfxsize=5in\centerline{\epsfbox{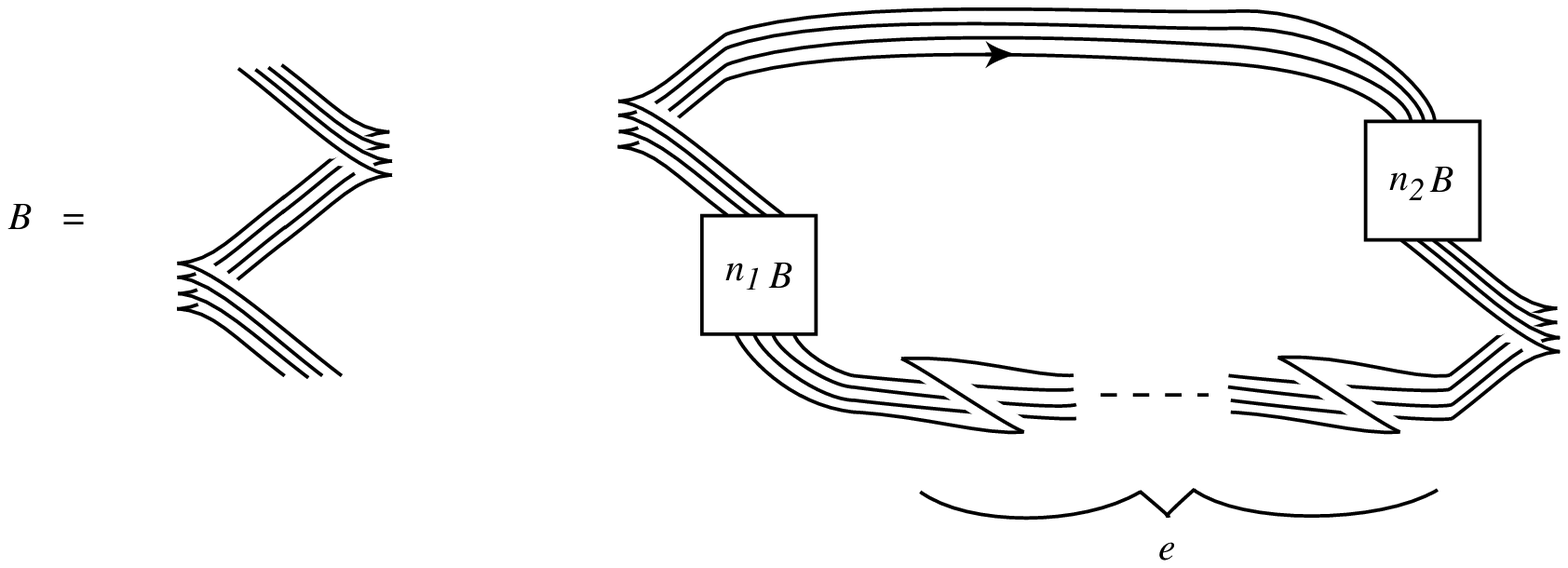}}}
	\caption{Legendrian torus knots.}
	\label{fig:torus}
\end{figure}

\begin{figure}[ht]
	{\epsfysize=1.4in\centerline{\epsfbox{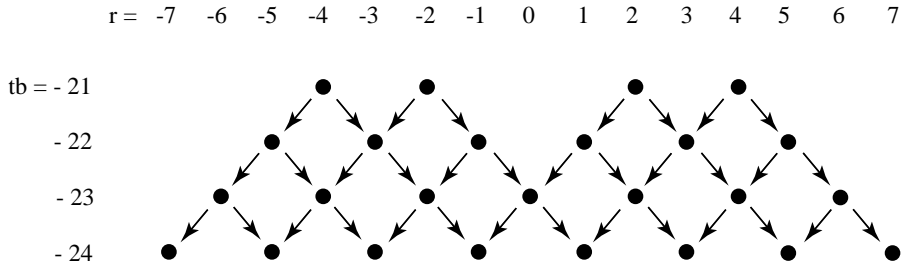}}}
	\caption{Some Thurston-Bennequin invariants and rotation numbers realized by Legendrian $(-7,3)$-torus knots.} 
	\label{fig:negexample}
\end{figure}

\begin{rem}
	{\em The above results imply that any Legendrian torus knot is some stabilization
	of a knot in Figure~\ref{fig:torus}. In Figure~\ref{fig:torus} note the positive torus knot 
	has $\tb=pq-p-q$ and $r=0$, while the negative $(p,q)$-torus knot has $\tb=pq$ and $r=q(n_2-n_1)-e,$
	where $|p|=nq+e$ and $n=1+n_1+n_2.$ To get an idea of the possible values of the invariants
	for negative torus knots consider the $(-7,3)$-torus knot. In Figure~\ref{fig:negexample} we show a graph of
	the (four largest) Thurston-Bennequin invariants versus the rotation numbers. Each dot corresponds
	to a pair of invariants that is realized by a Legendrian $(-7,3)$-torus knot and each arrow
	represents a stabilization.}
\end{rem}

In this section we classify positive, oriented Legendrian torus knots
in $(S^3,\xi)$ up to Legendrian isotopy, leaving the general case for the following subsection.
Precisely, in this section we prove:
\begin{thm}\label{thm:mainpositivelegendrian}
	If $K$ and $K'$ are two oriented Legendrian positive torus knots then
	they are Legendrian isotopic if and only if $\tb(K)=\tb(K')$ and
	$r(K)=r(K').$  Moreover, if $K$ is a $(p,q)$-torus knot, with $p,q>0$,
	then $\tb(K)\leq pq-p-q$ and if $\tb(K)=pq-p-q-n$
	for some non-negative integer $n$ then $r(K)\in\{-n,-n+2,\ldots,n\}.$
\end{thm}
We employ Strategy~\ref{mainstrat} described in Section~\ref{legendriansection} to prove this theorem.
First note that the maximal Thurston-Bennequin invariant for a Legendrian $(p,q)$-torus knot, with $p,q>0$, is
$\overline{\tb}_{(p,q)}=pq-p-q,$ since there is a Legendrian knot realizing this (see Figure~\ref{fig:torus}) and the
Bennequin inequality says that this is as large as it could be. Moreover, note that the Bennequin inequality
then implies that a Legendrian knot realizing this Thurston-Bennequin invariant must have rotation
number 0. With this in mind we show:
\begin{lem}\label{lem:uniquie+}
	Let $K$ and $K'$ be two topologically isotopic Legendrian positive torus knots with maximal 
	Thurston-Bennequin invariant. Then $K$ and $K'$ are Legendrian isotopic. 
\end{lem} 

\begin{proof}
Let $T$ and $T'$ be any standardly embedded tori in $S^3$ on which $K$ and $K'$ respectively sit.
We describe everything in terms of $T$ and $K$ but everything also holds for $T'$ and $K'.$
By Theorem~\ref{makeconvexrelcurve} we may make $T$ convex without moving $K$, since the
twisting of $K$ with respect to
$T$ is $-p-q$ (recall $p,q$ are positive). Now that $T$ is convex and $\tb(K)$ is maximal,
we may assume $T$ is in standard form.   This follows from observing that
$\#(K\cap \Gamma)> |K\cap \Gamma|$ implies the existence of a bypass along $K$, and hence a destabilization.
Here $\#(K\cap \Gamma)$ is the unsigned (actual) intersection number
and $|K\cap\Gamma|$ is the (minimum) geometric intersection number of the two isotopy classes.
Let $-\frac{r}{s},$ ($r,s>0$) be the slope of the dividing curves $\Gamma$ and $2n$ be the number of dividing curves.
According to Lemma~\ref{lem:invts}, $\tb(K)= pq-\frac{1}{2} \# (K\cap \Gamma)=pq-n\det\left( \begin{array}{cc}
p & -s\\ q & r\end{array} \right).$ Thus for $\tb(K)$ to equal $pq-p-q$ we must have $n=1$ and $r=s=1.$

Let $V'_0\cup V'_1$ and $V_0\cup V_1$ be the Heegaard splittings associated to
the tori $T'$ and $T$. Since the slopes of the dividing curves are the same, we may use the  classification
of tight contact structures on solid tori in Lemma~\ref{toriclassification} to find  
a contactomorphism $\phi:V_0\to V'_0$ (note that a meridional disk for $V_0$ of $V_0'$ has 
exactly one dividing curve). Applying Lemma~\ref{toriclassification}
again, we may extend $\phi$ to all of $S^3$, thus obtaining a contactomorphism of $S^3$ which takes $T$ to $T'.$
By Eliashberg's result 
(Theorem~\ref{thm:htpyequiv}) we may find a contact isotopy of $S^3$ taking $T$ to $T'.$
So now we may assume that $K$ and $K'$ are
two Legendrian knots on the same convex torus $T.$ Note $K$ and $K'$ are both leaves in the ruling foliation
of $T.$ The other ruling curves exhibit a Legendrian isotopy from $K$ to $K'.$
\end{proof}
The proof of Theorem~\ref{thm:mainpositivelegendrian} is thus complete by Lemma~\ref{aboutstabilization}
and:
\begin{lem}\label{lem:positivetorusstabilization}
	If $K$ is an oriented Legendrian $(p,q)$-torus knot, with $p,q>0$ and 
	$\tb(K)=pq-p-q-n,$ then there exist positive integers $n_1$ and $n_2$ such
	that $n=n_1+n_2,$ $r(K)=n_2-n_1$ and $K= S_-^{n_1}(S_+^{n_2}(K'))$ where $K'$
	is the unique Legendrian $(p,q)$-torus knot with maximal Thurston-Bennequin invariant. 
\end{lem}

\begin{proof}
Let $T$ be a standardly embedded torus on which $K$ sits. As in the proof of Lemma~\ref{lem:uniquie+}
we may assume that $T$ is convex and in standard form, since otherwise $K$ can be destabilized.
However, this time the slope $-\frac{r}{s}$ of the dividing curves
is not $-1$ or the number of dividing curves is not 2.
We first consider the case when $-\frac{r}{s}\not=-1.$ Let $V_0\cup V_1$ be the Heegaard splitting associated 
to the the torus $T.$ Recall the ``slope'' of the dividing curves is usually measured thinking of $T$ as
the boundary of $V_0.$ As the boundary of $V_0$ or $V_1$ the slope of the dividing curves on $T$ will be less than
$-1.$  Assume $V_0$ has this property. We know, by Lemma~\ref{solidtorustwisting}, that looking at concentric convex tori in 
$V_0$  we will see dividing curves with any slope in $[-\frac{r}{s},0).$ In particular, there will be
a torus $T'\subset V_0$ with two dividing curves having slope $-1.$ Let $U=T\times[0,1]$ be the region between
$T$ and $T'$ in $V_0$ and $A=K\times [0,1]$ be an annulus lying between $T=T\times\{0\}$ and $T'=T\times\{1\}$
in $U.$ The boundary of $U$ is convex and we may assume that the ruling curves on both boundary components
have slope $\frac{q}{p}.$ Thus we may assume that $\partial A=K\cup K'$ is Legendrian and $A$ is convex. 
The dividing curves will intersect $K,$ $N=2n\det\left( \begin{array}{cc}p & -s\\ q & r\end{array} \right)$ times 
and $K',$ $N'=2(p+q)$ times. As $r$ and $s$  are not both 1, $N'>N$, so we can find a boundary-parallel arc along
$K$ among the dividing curves of $A.$ This implies the existence of a bypass for $K$ and hence a destabilization.
Specifically, $K=S_\pm(K'')$ for some Legendrian $K''$. Repeating this argument we will eventually find
a sequence of destabilizations which will exhibit $K$ as $S_-^{n_1}(S_+^{n_2}(K'))$ where $K'$
is the unique Legendrian $(p,q)$-torus knot with maximal Thurston-Bennequin invariant. 

Now for the case when $-\frac{r}{s}=-1$ and $n>1.$ In this case we claim that there is a torus $T'$ in, say $V_1,$ parallel
to $T$ with two dividing curves having slope $-1.$ To see this just take a copy of $T$ inside $V_1$, make the ruling curves
meridional, and then look at a convex meridional disk $D.$ We may use the bypasses on $D$ to reduce the number of dividing
curves on the copy of $T$ until there are only two. (This follows easily from Proposition~\ref{prop:bypass}. Also see 
\cite{Honda1}.) 
Call the resulting torus $T'.$ Let $U=T\times [0,1]$ be the region
between $T$ and $T'.$ Now if we make the ruling curves on $T'$ have slope $\frac{q}{p}$ then using the annulus 
$A=K\times [0,1]$ in $U$ we may repeat the above argument to destabilize $K.$
\end{proof}

\subsection{Negative Legendrian Torus Knots}

We begin by establishing the upper bound on the Thurston-Bennequin invariant.

\begin{lem}
	For a nontrivial negative $(p,q)$-torus knot ($- p > q>0$) we have $\overline{\tb}_{(p,q)}=pq<0.$
\end{lem}

\begin{proof}
First note that the example in Figure~\ref{fig:torus} shows that $\overline{\tb}_{(p,q)}\geq pq.$
We show that $\overline{\tb}_{(p,q)}\leq pq$ by contradiction.
If $\overline{\tb}_{(p,q)}> pq$, then we can construct a Stein manifold $X$ (with boundary) by adding
a 2-handle to the 4-ball along a $(p,q)$-torus knot with framing $pq.$ (Given a symplectic 4-manifold with
convex boundary, one can add a symplectic 2-handle to a Legendrian knot $\gamma$ in the boundary with
framing $\tb(\gamma)-1$ and obtain a new symplectic 4-manifold with convex boundary, see \cite{el:fill, Gompf98}.)
According to \cite{Moser},  $\partial X$ is the connected sum of two lens spaces (neither of which is $S^3$
or $q=1$ and the knot is trivial). Now, a theorem of Eliashberg \cite{el:fill} says that $X$ must be a boundary sum
of two manifolds $X_0$ and $X_1.$ Using the Mayer-Vietoris sequence we see that one of the $X_i$'s, say $X_0,$ must be
an integral homology ball, but this is impossible since $\partial X_0$ is a nontrivial lens space.
\end{proof}

This establishes the second part of Theorem~\ref{maxtb}.  The easiest part of Strategy~\ref{mainstrat} to execute in
this situation is Part 2 concerning destabilization.

\begin{lem}
	If $K$ is an oriented Legendrian $(p,q)$-torus knot, with $pq<0$
	and $\tb(K)<pq,$ 
	then there is a $(p,q)$-torus knot $K'$ such that $\tb(K')>\tb(K)$ and $K$ is a stabilization
	of $K'.$
\end{lem}

\begin{proof}
Since $\tb(K)\leq pq$ we can use Theorem~\ref{makeconvexrelcurve} to show $K$ lies on a convex standardly embedded torus $T.$
We know the dividing curves $\Gamma$ on $T$ have slope $-\frac{r}{s}\not=\frac{q}{p},$ since
\begin{equation}
	\tb(K)=pq-\frac{1}{2}\#(\Gamma\cap K)<pq.
\end{equation}
Now as measured on either $V_0$ or $V_1$ the dividing curve have slope less than $\frac{q}{p}.$ Assume
$V_0$ has such dividing curves. Now since Lemma~\ref{solidtorustwisting} says that we can find tori 
in $V_0$ whose dividing curves have any slope in  
$[-\frac{r}{s}, 0)$ we can find a torus $T'$ in $V_0$ whose dividing curves have slope $\frac{q}{p}.$ Now as in
the proof of Lemma~\ref{lem:positivetorusstabilization} we can take an annulus $A$ between $T$ and $T'$ with one 
boundary on $K$ and the other on a Legendrian divide of $T'$ and thus find a bypass for $K$ since the dividing curves
on $A$ will intersect the boundary component containing $K$ more than the one containing the Legendrian divide.
Hence we may use this bypass to destabilize $K.$   
\end{proof}

Before proceeding to Parts 1 and 3 of Strategy~\ref{mainstrat} we need to take a detour concerning
rotation numbers and relative ``Euler classes.''  In the following we will assume that $K$ is a Legendrian
$(p,q)$-torus knot with maximal Thurston-Bennequin invariant.
Let $T$ be a convex standardly embedded torus on which
$K$ sits (Theorem~\ref{makeconvexrelcurve} assures we can find such a torus). We may put $T_\xi$ in standard
form with $K$ being one of the Legendrian divides.
We note that all the Legendrian divides of $T$ are Legendrian isotopic. This follows from Lemma~\ref{model},
since all the Legendrian divides sit on a linearly foliated torus. 
Using an argument like the one at the end of the proof of Lemma~\ref{lem:positivetorusstabilization},
we may now reduce the number of Legendrian divides on $T$ to two and assume that $K$ is one of these divides.

Now let $S^3=V_0\cup_T V_1$, where $V_0$ (resp. $V_1$)
is the solid torus with meridional curve $\mu$ (resp. $\lambda$).
From Lemma~\ref{solidtorustwisting} we know that inside $V_0$ there is a solid torus $S$ with two dividing curves
of slope $-\frac{1}{m+1},$  where $|p|=mq+e,$ and there is a solid torus $S'$ containing $V_0$ with two
dividing curves of slope $-\frac{1}{m}.$ Let $T_m=\partial S'$ and $T_{m+1}=\partial S.$ In addition, set
$\overline{S}=\overline{S^3\setminus S'}$ and $\overline{S'}=\overline{S^3\setminus S}.$

Now we define an invariant of homology classes of curves on $T.$ Let $v$ be any globally non-zero section
of $\xi$ and $w$ a section of $\xi\vert_T$ that is transverse to and twists (with $\xi$) along the 
Legendrian ruling curves and is tangent to the Legendrian divides. If $\gamma$ is a closed oriented curve on
$T$ then set $f_T(\gamma)$ equal to the rotation of $v$ relative $w$ along $\gamma.$ One may check the following 
properties (cf. \cite{Etnyre99}).
\begin{itemize}
	\item The function $f_T$ is  well-defined on homology classes.
	\item The function $f_T$ is unchanged if we isotop $T$ among convex tori in standard form.
	\item If $\gamma$ is a $(r,s)$-ruling curve or Legendrian divide then $f_T(\gamma)=r(\gamma).$ 
\end{itemize}
We may similarly define $f_m$ and $f_{m+1}$ for curves on $T_m$ and $T_{m+1}.$ 
The main facts we need concerning these invariants are:
\begin{enumerate}
	\item $f_T(\mu) = 1-q$ or $q-1.$
	\item $f_m(\lambda)$ is in $\{m-1,m-3,\ldots, 1-m\}.$
	\item If $f_T(\mu)=1-q$ then $f_T(\lambda)=f_m(\lambda) + (m-|p|).$ So $f_T(\lambda)$ is in $\{2m-|p|-1,
		2m-|p|-3,\ldots,
		1-|p|\}.$
	\item If $f_T(\mu)=q-1$ then $f_T(\lambda)=f_m(\lambda) + (|p|-m).$ So $f_T(\lambda)$ is in $\{|p|-1,|p|-3,\ldots,
		|p|-2m+1\}.$

\end{enumerate}
We will prove these facts at the end of this section.
From the above properties we know that 
\begin{equation}\label{rfromf}
	r(K)=pf(\mu)+qf(\lambda).
\end{equation}
Thus the possible values for $r(K)$ lie in 
$$\{ \pm(\vert p\vert-\vert q\vert - 2qk) : k\in\mathbb{Z},
0\leq k<  \frac{\vert p\vert-\vert q\vert}{\vert q\vert}\}.$$
Since all these possible rotation numbers
are realized by Legendrian knots in Figure~\ref{fig:torus}, this finishes the proof of Theorem~\ref{rrange}.

\begin{lem}\label{classifynegtoursknots}
	Let $K$ and $K'$ be two topologically isotopic Legendrian negative torus knots with maximal 
	Thurston-Bennequin invariant. Then $K$ and $K'$ are Legendrian isotopic if and only if $r(K)=r(K')$. 
\end{lem}

\begin{proof}
Let $T$ and $T'$ be tori on which $K$ and $K'$ respectively sit. 
We assume that $T$ and $T'$ have been arranged as described above.
The proof is finished as we finished the proof of Lemma~\ref{lem:uniquie+}. We only need to recall  that
Theorem~\ref{toriclassification}
says that the contactomorphism type of a tight contact structure on a solid torus (with standard convex boundary) 
is determined by the number of positive bypasses on a meridional disk
and the number of positive bypasses on
the meridional disks for $V_0$ and $V_1,$ where $S^3=V_0\cup_T V_1,$ are determined by the rotation number of $K.$
\end{proof}

We also may show:
\begin{lem}\label{valleys}
	Recall $m$ and $e$ are integers such that $|p|=mq+e$ where $0<e<q.$
	Let $K$  and $K'$ be negative Legendrian $(p,q)$-torus knots with maximal Thurston-Bennequin
	invariant. 
	If the rotation numbers are $r$ and $r-2e,$ respectively, 
	then $S^e_-(K)$ and $S^e_+(K')$ are Legendrian isotopic. 
	If the rotation numbers are  $r$ and $r-2(q-e),$ respectively, 
	then $S^{(q-e)}_-(K)$ and $S^{(q-e)}_+(K')$  are Legendrian isotopic.
\end{lem}

\begin{proof}
We assume $r(K)=r(K')+2e$ and leave the similar case to the reader. In this situation there must be a
$k$ so that $f_T(\lambda)=k+(m-|p|)$ and $f_{T'}(\lambda)=k+(|p|-m).$ Thus $f_T(\mu)=1-q$ and $f_{T'}(\mu)
=q-1.$ In the notation set up above, let $T_m$ we the convex torus outside $V_0$ with boundary slope $-\frac{1}{m}$ and
$T_m'$ be the corresponding one for $T'.$ From Facts  3 and 4 stated above we know that $f_{T_m}(\lambda)=k$ and
$f_{T'_m}(\lambda)=k.$ By examining dividing curves on the meridional disks, we also have $f_{T_m}(\mu)=f_{T'_m}(\mu)=0$
(from Equation~\ref{formulafore}).
Now arrange for the Legendrian ruling curves on $T_m$ to have slope $\frac{q}{p}$ and  consider the annulus $A$
between $T$ and $T_m$ with slope $\frac{q}{p}.$ We can take one boundary of $A$ to be $K\subset T$ and the
other a Legendrian ruling curve on $T_m.$ Once we make $A$ convex, the dividing curves do not intersect the boundary
component touching
$T$ and intersect the other boundary component $2\det\left( \begin{array}{cc} p & -m\\ q & 1\end{array} \right)
=2(p+qm)=-2e$ times. Since the rotation number (values of $f_T$ and $f_{T_m}$) of the boundary components of $A$
differ by $-e$, we may use Lemma~\ref{lem:invts} (or a slight generalization of it, see \cite{EH}) to see there are
$e$ boundary-parallel dividing curves separating off $e$ negative disks from $A.$ This shows that
$K_s=A\cap T_m=S_-^e(K).$ Similarly, we see that $K'_s\subset T'_m$ is $K_+^e(K').$ We may now use the argument in
the proof of Lemma~\ref{classifynegtoursknots} to find a contact isotopy taking $T_m$ to $T_m'.$  A further contact
isotopy takes $K_s=S_-^e(K)$ to $K'_s=K_+^e(K').$ 
\end{proof}

\begin{thm}
	Negative Legendrian torus knots are determined by their knot type, Thurston-Bennequin invariant
	and rotation number.
\end{thm}
With this theorem we have completed the proof of Theorem~\ref{main} and the classification of Legendrian torus knots.
Before we begin the proof note that if we graph the Thurston-Bennequin invariants versus the rotation numbers 
that a fixed negative torus knot
realizes, as in Figure~\ref{fig:negexample}, we obtain a figure that looks like a mountain range. Above we
have shown that the dots at the peak of the mountains (corresponding to maximal Thurston-Bennequin invariant
knots) are unique and the dots in the valleys (where stabilizations of maximal Thurstion-Bennequin invariant 
knots first have the same invariants, in Figure~\ref{fig:negexample} there are the points $(-22,-3), (-23,0)$ and $(-22,3)$) 
are unique. The proof of this theorem follows from this and the properties of stabilizations.
\begin{proof}
Let $K$ and $K'$ be two $(p,q)$-torus knots with the same invariants. Let $K$ destabilize to $K_d$ and $K'$
to $K_d'.$ If $K_d$ and $K_d'$ have the same invariants, then they are Legendrian isotopic and $K$ and $K'$ are
the same stabilization of the same knot and hence are Legendrian isotopic. Now suppose the rotation numbers
of $K_d$ and $K_d'$ differ by $2e$ as in Lemma~\ref{valleys}. Then we can realize $K$ as a stabilization of $S_-^e(K_d)$ 
and $K'$
as a stabilization of $S_+^e(K_d').$ Thus they both destabilize to $S_-^e(K_d)=S_+^e(K_d')$ and hence they {\em{both}}
destabilize to $K_d$ and $K_d',$ implying they are Legendrian isotopic. Finally, if the rotation numbers of
$K_d$ and $K_d'$ differ by something else, then a similar argument will show that $K$ and $K'$ will destabilize
to the same Legendrian knot, finishing the proof.
\end{proof}

\begin{proof}[Proof of Fact 1:]
Arrange for all the Legendrian ruling curves on $T_m,$ 
$T_{m+1}$, and $T$ to be meridional and let $D$ be a meridional disk for $S'$ that intersects these
three tori in Legendrian curves. We may now isotop $D$ (relative to its intersection with the tori)
so that it is convex. Let $D=D_V\cup A$ where $D_V=D\cap V_0$ and $A=\overline{D\setminus D_V}.$ 
Orient $D$ so that $\mu$ is the oriented boundary of $D_V.$ Now the dividing curves on $D_V$ intersect
$\mu=\partial D_V,$ $2q$ times (since that is the number of times the dividing curves on $T$ intersect $\mu$).
Moreover, they intersect each of $\partial D$ and $D\cap T_m,$ 2 times. 
We claim that all the dividing curves on $D_V$ separate off disks that contain no dividing curves.
Note this implies that all the bypasses on $D_V$ are of the same sign. If this were not the case, then
we would have bypasses on $D_V$ of both signs and we would be able to glue this to a bypass of the same
sign on $A$,  creating an overtwisted disk. To match a bypass on $D_V$ with any bypass on $A$ we might have to ``add copies of
$T$ to $V_D$,'' that is, take a copy of $T$, cut it along $\mu$ to obtain an annulus, and glue one of its boundary
components to $V_D$ and the other to $A$. This has the effect of shifting the dividing curves on $D_V$ relative
to those on $A$ by $\frac{q}{p}$ (for details see the section on Sliding Maneuvers in \cite{Honda1}).
Since all the bypasses have the
same sign, Equation~\ref{formulaforrfromconvexity} implies that $f_T(\mu)=r(\mu)=-1+q$ or $1-q.$
\end{proof}

\begin{proof}[Proof of Fact 2:]
Make the Legendrian ruling curves on $T_m=\partial\overline{S}$ be meridional and
let $D$ be the meridional disk for $\overline{S}$ with Legendrian boundary.
Note $\partial D=\lambda$ (or at least a translate of it). Now the dividing curves 
intersect the boundary of $D,$ $2m$ times. And since there are no closed homotopically trivial dividing curves 
we can conclude that there are exactly $m$ dividing curves.
By examining the  possible configurations and using Equation~\ref{formulaforrfromconvexity}, one may readily
conclude that $f_m(\lambda)\in\{m-1,m-3,\ldots, 1-m\}.$ Moreover, by recalling how $f_m(\lambda)$ is related
to the rotation number of $K$ and noting the values of $r(K)$ realized in Figure~\ref{fig:torus}, we see that
all the possible values of $f_m(\lambda)$ are actually realized.
\end{proof}

\begin{proof}[Proof of Fact 3:]
Make the Legendrian ruling curves on $T_m$ and $T_{m+1}$ longitudinal and 
let $A$ be a longitudinal annulus spanning between $T_m$ and $T_{m+1}$ with Legendrian boundary.
After making $A$ convex, the dividing curves will intersect $S^1_m=T_m\cap A$ in $2m$ points and
$S^1_{m+1}=T_{m+1}\cap A$ in $2m+2$ points. We claim that $2m$ dividing curves run from one boundary
component of $A$ to the other boundary component and there is one boundary-parallel dividing curve
with endpoints on $S^1_{m+1}.$ If this is not the case, then there must be boundary-parallel dividing
curves for $S^1_m$ and hence a bypasses. Carefully applying Lemma~\ref{bypassalteration}, we can then find
a convex torus between $T_m$ and $T_{m+1}$ with slope not lying between $[-\frac{1}{m},-\frac{1}{m+1}]$.
This, by Lemma~\ref{torustwisting}, implies that there is a convex torus with boundary slope 0 and the
Legendrian divides on this torus will be the boundaries of overtwisted disks. Therefore, the dividing curve
configuration on $A$ is as claimed.

Now make the Legendrian ruling curves on $T$ longitudinal and make $A$ intersect $T$ in one of these
longitudinal curves, say $\gamma.$ The curve $\gamma$ separates $A$ into two annuli --- $A_m$ which
touches $T_m$ and $A_{m+1}$ which touches $T_{m+1}.$ Note the dividing curves are still as described
above. From this we can deduce the structure of the dividing curves on $A_m.$ To this end, note
the dividing curves on $A_m$ intersect the boundary component touching $T_m,$ $2m$ times and 
the boundary component touching $T,$ $2|p|$ times. Due to the structure of the dividing curves on $A$,
we know that the $2m$ dividing curves emanating from $T_m$ run across $A_m$ to the other boundary 
component. Thus we know there are $|p|-m$ other dividing curves, all of whose boundaries lie on $T.$
As argued in the proof of Fact 1, we can conclude that these all separate off disks that contain no dividing curves
and thus give $|p|-m$ bypasses of the same sign.

Let $e_U$ be the relative Euler class for the region $U=T^2\times I$ between $T$ and $T_m$.
Since $e_U(\lambda)=f_T(\lambda)-f_m(\lambda)$, Equation~\ref{formulafore} yields
$$f_T(\lambda)-f_m(\lambda)= (m-|p|) \hbox{ or } (|p|-m).$$
If $f_T(\lambda)=1-q$, then we claim that $f_T(\lambda)-f_m(\lambda)=m-|p|.$  Assuming the contrary,
we obtain a contradiction to the tightness of $\xi$ as follows: 
Using Equation~\ref{formulaforrfromconvexity} one may easily check that $f_m(\mu)=0.$ Thus,
if $f_T(\mu)=1-q=f_T(\mu)-f_m(\mu)=e_U(\mu),$ then $e_U(p\mu+q\lambda)=2q(|p|-m)-e>e.$ However, looking
at the convex annulus $A'$ of slope $\frac{p}{q}$ in $U$, we see that the there are only $e$
non-closed dividing curves, all of which have their boundary on $T_m.$ This means that for $e_U(\mu)>e$ there
must be some dividing curves bounding disks, violating the tightness of $\xi.$ Therefore, $f_T(\lambda)-f_m(\lambda)=m-|p|.$
\end{proof}

\begin{proof}[Proof of Fact 4:]
This is quite similar to the proof of Fact 3 and is left to the reader.
\end{proof}

\subsection{Transversal Torus Knots}

Since Theorem~\ref{main} implies that Legendrian torus knots are stably simple, we may use
Theorem~\ref{thm:stably-transversally} to conclude
\begin{cor}
	If $K$ and $K'$ are transversal torus knots then
	there are transversally isotopic if and only if $l(K)=l(K').$ 
	If $K$ is a $(p,q)$-torus knot, with $p,q>0,$ then
	$l(K)\leq pq-p-q$ and if $p<0$, $q>0$ then $l(K)=pq+\vert p\vert-\vert q\vert.$
	Moreover, all odd integers satisfying these bounds are realized by 
	the appropriate transversal torus knots.
\end{cor}
As mentioned earlier, this extends \cite{Etnyre99} and reproves a special case of a result of
Birman and Wrinkle \cite{bw} and Menasco \cite{M}.

\section{Figure Eight Knots}

The following is the main theorem of this section:

\begin{thm}  \label{thm:maineight} The figure eight knot type $\mathcal{K}$ is Legendrian simple.
\end{thm}

Again, Strategy \ref{mainstrat} will be employed to prove Theorem \ref{thm:maineight}.   


\subsection{Determination of $\overline{\tb}$}

The first step in the proof of Theorem \ref{thm:maineight} is to compute the maximal
Thurston-Bennequin invariant $\overline{\tb}$ for the figure eight knot.

\begin{lem}\label{lem:maxtb}
The maximal Thurston-Bennequin invariant for the figure eight knot is $-3$.
\end{lem}

\begin{figure}
	{\epsfysize=1in\centerline{\epsfbox{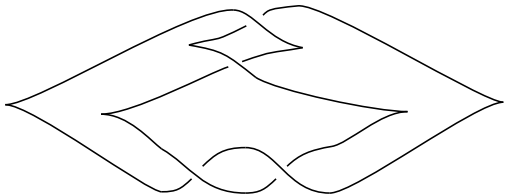}}}
	\caption{Legendrian figure eight knot with $\tb=-3$.}
	\label{legfigeight}
\end{figure}

\proof
To prove the upper bound $\overline{\tb}\leq -3$, use
the inequalities due to Fuchs-Tabachnikov \cite{ft,Ta}.  They give upper bounds for
$\tb+r$ of a Legendrian knot in terms of the HOMFLY and Kauffman polynomials.  A simple calculation of
these polynomials give the desired upper bound.  Figure \ref{legfigeight}
provides an explicit example of
a Legendrian figure eight knot with $\tb=-3$.
\qed

\subsection{Unique maximal $\tb$ representative}

In this section we prove the following theorem:

\begin{thm}\label{thm:uniquemax}	
	If $K$ and $K'$ are two oriented Legendrian figure eight $\mbox{knots}$ with
	maximal Thurston-Bennequin invariant  $\tb=-3$
	{\em (}and necessarily $r(K)$ $=r(K')=0${\em )}, then $K$ and $K'$ are Legendrian isotopic.
\end{thm}

Let $M$ be $S^3\setminus N$
where $N$ is a small standard neighborhood of $K$ so that the Legendrian ruling on $\partial M$
is by longitudes (this is possible since $\tb<0$)  --- note that there is a well-defined longitude
coming from the Seifert surface.   Recall we have a fibration $p\co M\to S^1$
with fiber a punctured torus $\Sigma$ with Legendrian boundary.

The following is a general procedure for proving that two Legendrian knots $K$ and $K'$
of the same knot type
and maximal Thurston-Bennequin invariant are isotopic.
We construct a {\it contactomorphism}
$f\co S^3 \to S^3$ taking $K$ to $K'$ as follows:  Let $f$ be a map from a standard neighborhood $N$
(as above) of $K$ to a standard neighborhood $N'$ of $K'$ --- this is possible because of local
uniqueness. Assume $N$ and $N'$ have already been identified via $f$.  Then we have two tight
contact structures on $M$ with identical boundary characteristic foliations.
Assume we can prove the following:
\begin{claim}
The two tight contact structures on $M$ are contactomorphic relative the boundary.
\end{claim}

Then the map $f$ may be extended to a contactomorphism of $S^3$ onto $S^3$, taking $K$ to $K'$.
With the map $f$ in hand we may use Theorem \ref{thm:htpyequiv} to show we have a contact isotopy
from the identity map $id: S^3\rightarrow S^3$ to $f$, exhibiting the desired Legendrian
isotopy from $K$ to $K'.$

We have therefore reduced the proof of Theorem \ref{thm:uniquemax} to the analysis (and normalization)
of the tight contact
structure on $M=S^3\setminus N$.  
Since $M$ fibers over the circle with  punctured
torus fiber $\Sigma=T^2\backslash D^2$, we write $M=(\Sigma \times I)/\sim$, where
$(x,0)\sim (\Psi(x),1)$, and
$$\Psi=\begin{bmatrix} 2 & 1 \\ 1 & 1 \end{bmatrix}$$ is the monodromy map.  
$\Psi$ is a diffeomorphism of $\Sigma$ which fixes $\partial D^2$.  Note that the matrix by itself
does not completely define $\Psi$, since it is possible to compose with Dehn twists along curves 
parallel to $\partial D^2$. Since 
$$\begin{bmatrix} 2 & 1 \\ 1 & 1 \end{bmatrix}=\begin{bmatrix} 1 & 1 \\ 0 & 1 \end{bmatrix}
\begin{bmatrix} 1 & 0 \\ 1 & 1 \end{bmatrix},$$
we define $\Psi$ as a Dehn twist along a closed curve parallel to $(0,1)$, followed by a 
Dehn twist along a closed curve parallel to $(1,0)$. In particular, we compose with no Dehn twists along $\bdry D^2$.
We will use this particular matrix representative of the conjugacy class of $\Psi$.

Consider the characteristic foliation on a convex $\Sigma$ with Legendrian boundary.
The following several pages will be devoted to  normalizing $\Sigma$.  Note that the dividing set
$\Gamma_\Sigma$ will consist of closed curves or arcs with endpoints on $\bdry \Sigma$.

First observe the following simple lemma:
\begin{lem}\label{lem:curves}
	If $\{\gamma_1,\ldots\gamma_n\}$ is a collection of non-boundary-parallel, mutually
	non-parallel arcs  in a punctured torus $\Sigma$ (with endpoints on $\bdry \Sigma$),
	then $n\leq 3.$  Moreover, if there exists a closed curve on $\Sigma$ disjoint from
	the $\gamma_i$'s, then $n\leq 1.$
\end{lem}

The possible configurations of dividing curves are given below:

\begin{prop}\label{prop:curves}
	Suppose the dividing curves on a convex $\Sigma\subset M$ with Legendrian
	boundary consist of arcs
	$\gamma_1,\ldots,\gamma_n$
	with endpoints on $\bdry \Sigma$, and closed curves $c_1,\ldots,c_m.$
	\begin{itemize}
		\item[I.] If the $\gamma_i$'s realize one distinct non-boundary-parallel isotopy class
			on $\Sigma,$ then $\Sigma$ may be isotoped so that  {\em(a)} $n\geq 1$ {\em(}$n$ odd{\em)} and $m=1$,
			{\em (b)} $n=0$ and $m=2$, or {\em(c)} $\Sigma$ has a bypass along
			$\bdry \Sigma$.
		\item[II.] If the $\gamma_i$'s realize two distinct non-boundary-parallel isotopy classes
			on $\Sigma,$
			then $m=0$ and either {\em(a)} each  isotopy class has two $\gamma_i$, or {\em(b)}
			$\Sigma$ may be isotoped so that it has a bypass along
			$\bdry \Sigma$.
		\item[III.] If the $\gamma_i$'s realize three distinct non-boundary-parallel isotopy classes
			on $\Sigma,$
			then $m=0$ and either {\em(a)} one of the isotopy classes has one $\gamma_i$ and the other two
			have one or more, or {\em(b)} $\Sigma$ may be isotoped so
			that it has a bypass along $\bdry \Sigma$.	
	\end{itemize}
\end{prop}

The proposition is preliminary.  It will be improved later, after considerably more work.

\begin{proof}
Note that any claim we will be making concerning $\Sigma$ should
be interpreted up to a $C^0$-small isotopy --- in particular, we will make extensive use of the
Legendrian Realization Principle (Theorem \ref{lerp}) without explicitly mentioning each time that 
a $C^0$-small perturbation is taking place first.  
Let $\Sigma$ be a convex fiber in $M$ with Legendrian boundary, and $N$ an open convex
neighborhood of $\Sigma$.
Then in $(M\setminus N)=\Sigma\times [0,1],$
we look for an annulus $A$
with $\partial A=A_0- A_1$, where the $A_i$ are Legendrian curves in (a $C^0$-isotoped
copy of) $\Sigma_i=\Sigma\times\{i\}$, and $A$ contains a bypass.  Denote the dividing set of $\Sigma_i$
by $\Gamma_i$.

Let us consider Case I.  Observe that $m+n$ is always even.
Here we may always find a closed Legendrian curve $c$ on
$\Sigma_0$ which does not intersect $\Gamma_0$.
We take $A$ to be the annulus $c\times [0,1].$  Note that $\Gamma_1=\Psi(\Gamma_0)$
and the geometric intersection number
$|(c\times\{1\}) \cap \Gamma_1| \geq m+n$, because an Anosov map cannot fix a
curve isotopy class.   Realize $c\times\{1\}$ as a Legendrian curve which has minimal
intersection number with $\Gamma_1$ in its isotopy class.
Then $A$ will have at least $n+m$ dividing curves emanating from $A_1$ and none from $A_0$, from which we can
conclude using the Imbalance Principle (Proposition \ref{imba}) that there is a bypass along $A_1.$

If $m+n\geq 4$, then we may always find an
embedded bypass, and we apply Proposition \ref{prop:bypass}.  Let
$\delta$ be the attaching Legendrian arc $\subset \Sigma_1$ for the bypass, and $p,q,r$ be the
intersections with $\Gamma_1$ in consecutive order along $\delta$.  
Since $m+n\geq 4$, $p,q,r$ lie on distinct dividing curves.
If $p,q$ both lie on the $\gamma_i$'s  (or the same for $q,r$),  then there will exist a boundary-parallel
dividing curve on the new $\Sigma_1$, after bypass attachment (and thus a bypass for the new $\partial \Sigma_1$).
If $p,q$ both lie on the $c_i$'s (or the same for $q,r$), then  we will be able to reduce $m$
by two.   Hence, if $m\geq 2$ and $n\geq 2$, every bypass will be as above, and we may either reduce $m$ by 2 or 
obtain a bypass for $\Sigma_1$ along $\bdry \Sigma_1$.
Now there are two other cases:  $n\geq 3$, $m=1$, and $p,r$ lie on $\gamma_i$'s and $q$ lies on a $c_i$;
or $n=1$, $m\geq 3$, and $p,r$ lie on $c_i$'s and $q$ lies on a $\gamma_i$.  In the latter case, we reduce $m$ by 2.
(Note that, in the former case, the bypass attachment yields a dividing set which puts us in Case III, where
two of the isotopy classes of arcs are represented by a single $\gamma_i$.)

We eventually have (a) $n\geq 1$ and $m=1$, (b) $n=0$ and
$m=2$, (b') $n=2$ and
$m=0$, or (c) a bypass along $\bdry \Sigma_1$.
We can actually do better, and remove (b') $n=2$, $m=0$.  If $A_1$ intersects $\Gamma_1$ exactly twice,
then we have a degenerate bypass along $A_1$ (one with same endpoints), and attaching this degenerate
bypass gives rise to a dividing set which has a boundary-parallel dividing curve.
If the intersection number is greater than 2, then keep $A$ but repeatedly attach bypasses until we
arrive at boundary-parallel components of (the new) $\Gamma_1$.
Observe that, once we have reduced to (a) or (b), there will still exist bypasses that can be attached onto $\Sigma_1$.
However, they do not necessarily yield
boundary-parallel dividing curves after attachment and
cannot be used to destabilize $\partial \Sigma_1$ --- at least not immediately.

In Cases II and III,  $m=0$ by Lemma \ref{lem:curves}.   Consider Case II.
Let $\gamma_1$ and $\gamma_2$ be the two isotopy classes of curves of $\Gamma_0$, and $n_1$, $n_2$ be the number
of arcs in each. Note that both $n_1$ and $n_2$ are even. 
For convenience, we will identify $\Sigma=(\R^2/\Z^2)\backslash D^2$, and $\gamma_1$, $\gamma_2$ with
minimal (shortest-length) integral vectors $v_1, v_2\in \Z^2$.  
It is easy to find a curve $c$ on $\Sigma$ so that $c$ intersects
$\Gamma_0$ and $\Gamma_1$ in a different number of points.
Indeed consider the subset $\mathcal{C}_{min}$ of the set $\mathcal{C}$ of (isotopy classes
of) closed (connected) curves $c$ on $\Sigma$,
consisting of $c$ which
have minimal geometric intersection $|c\cap \Gamma_0|$ among curves in $\mathcal{C}$.
Since the Anosov map $\Psi$ cannot preserve a finite set
of closed curves, $\mathcal{C}_{min}\not=\Psi(\mathcal{C}_{min})$.  Noting that
$\Psi(\mathcal{C}_{min})$ consists of curves $c$ with minimal intersection number 
$|c\cap \Gamma_1|$ among curves in $\mathcal{C}$, we pick $c\in \mathcal{C}_{min}\backslash
\Psi(\mathcal{C}_{min})$.
Using this $c$ we may argue as in Case I to find a bypass along $\Sigma_1$.

As long as $n_1, n_2>2$, $p,q,r$ (as above) lie on three distinct dividing curves, and $p,q$
(or $q,r$) lie on `consecutive' dividing curves of the same type $v_1$ (or $v_2$).
Thus we can produce a boundary-parallel component of the new $\Gamma_1$ after attaching the bypass to
$\Sigma_1$. Therefore, at least one of the $n_i$ must equal $2$.
We will eliminate $n_1=2$, $n_2>2$ (or $n_1>2$, $n_2=2$).  In this case,
$\mathcal{C}_{min}$ consists of one curve $d_2$ parallel to $v_2$.  As above, $d_2\not \in \Psi
(\mathcal{C}_{min})$, so $|\Gamma_1\cap d_2|\geq n_2\geq 4$, since $d_2$ intersects both
$\Psi(v_1)$ and $\Psi(v_2)$.
We then let  $A=d_2\times[0,1]$, and find a bypass along $d_2\times \{1\}$.
The only nontrivial case is when our bypass along $\Sigma_1$ involves
only the two dividing curves isotopic to $\Psi(v_1)$,  $p,r$ lie on the same curve, and the bypass attachment
changes the dividing curves $\Psi(v_1)$, $\Psi(v_2)$  to $\Psi(v_1\pm v_2)$, $\Psi(v_2)$.
In this case we will have a sequence of `nested bypasses'.  More precisely,
attach the (outermost) bypass onto $\Sigma_1$ and take a new $A=d_2\times[0,1]$ such that
$d_2$ has fewer intersections with the new $\Gamma_1$.  Note that $d_2\times \{0\}$ intersects
$\Gamma_0$ twice, but $d_2\times \{1\}$ still intersects $\Gamma_1$ at least four times.
In any event, if we reduce this intersection number to four, there can no longer exist any bypasses of
the nontrivial type, and we produce a destabilization.
The only case left is when $n_1=n_2=2$.

In Case III, if we crush the
boundary $\Sigma$ to a point, then homologically one of the isotopy classes of $\Gamma$ is the sum of the
other two isotopy classes.  
Therefore we write the three classes as $v_1$, $v_2$, $v_1+v_2$, and the
number of dividing curves as $n_1$, $n_2$, $n_{12}$
(these must have the same
parity).   Using the same argument as in Case II (with $\mathcal{C}_{min}$), we produce an annulus $A$, together with a bypass
along $\Sigma_1$.

As long as three of the isotopy classes have more than one dividing curve,
any bypass along $A_1$ gives rise to a boundary-parallel component of $\Gamma_1$ on $\Sigma_1$.
Thus, one of the isotopy classes will have one dividing curve.
\end{proof}

Returning now to our normalization, we know that since $\tb(K)$
is maximal, there can be no bypasses on $\Sigma$. Thus Proposition \ref{prop:curves} and
Lemma \ref{lem:maxtb} imply that there are precisely three dividing curves --- we have either
I(a) with $n=3$, $m=1$ or III(a).  In what follows, we need to deal with configurations of both types,
but, for the time being, we may assume that all of the dividing curves are arcs and
they are non-parallel (we can switch from I(a) to III(a) according to the proof of the previous proposition).
Let us still denote the three arcs by the corresponding minimal integral basis vectors $v_1$, $v_2$, 
together with their sum $v_1+v_2$. 

\begin{figure}
	{\epsfysize=2in\centerline{\epsfbox{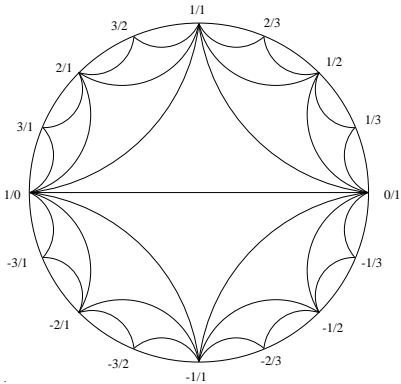}}}
	\caption{The standard tessellation of the hyperbolic unit disk}
	\label{tessellation}
\end{figure}

We will now normalize these vectors, showing that we can always assume we are in one of two possible
situations. To do this we need some way of listing all possible triples of such vectors (or, equivalently,
all possible two-triangle ``triangulations'' of $T^2$). We do this using the standard (Farey) tessellation
of the hyperbolic unit disk $D^2=
\{(x,y)|x^2+y^2\leq 1\}$, which we now review.
Recall we start by labeling $(1,0)$ as $0={0\over 1}$, and $(-1,0)$ as $\infty={1\over 0}$. 
We inductively label points on $S^1=\bdry D^2$ as follows (for $y>0$):  Suppose we have already
labeled $\infty \geq {p\over q}\geq 0$ ($p,q$ relatively prime) and $\infty \geq {p'\over q'}\geq 0$ ($p',q'$ relatively prime)
such that $(p,q)$, $(p',q')$ form a $\Z$-basis
of $\Z^2$.  Then, halfway between ${p\over q}$ and ${p'\over q'}$ along $S^1$ on the shorter arc (the one for which
$y>0$), we label ${p+p'\over q+q'}$.
We then connect two points ${p\over q}$ and ${p'\over q'}$ on the boundary by a hyperbolic geodesic, 
if the corresponding shortest integral
vectors form an integral basis of $\Z^2$.  See Figure \ref{tessellation}.
An important observation for us is that the tessellation gives a triangulation of $D^2$ 
into a union of ideal triangles, and that each
geodesic is on the boundary of exactly two ideal 
triangles --- any two integral basis vectors $v_1$, $v_2$ can be completed into
a triple
$\{v_1,v_2,v_3\}$ which mutually form an integral basis, in exactly two ways.
Observe that the three vectors $v_1, v_2, v_1+v_2$ in the previous paragraph form the vertices of an ideal triangle
in the standard tessellation and any ideal triangle gives a possible triple.
Note also that the vertices  of the triangles in the tessellation are labeled with the slopes of the vectors in the
corresponding triangulation of $T^2.$ More generally, there is a correspondence between the points on the 
``circle at infinity'' $S^1_\infty$ of the hyperbolic disk and slopes on $T^2.$

\begin{prop} \label{normalization}
There exists an isotopic copy of $\Sigma$ with dividing set $\Gamma=$ $\{(1,1),$ $(0,1),$ $(1,2)\}$ or
$\{(1,0),(0,1),(1,1)\}$.
\end{prop}

To prove proposition~\ref{normalization} we start with $\Gamma_0$ consisting of three
dividing curves corresponding to $v_1$, $v_2$, $v_1+v_2$, where $\{v_1,v_2\}$ form an integral 
basis of $\Z^2$ and thus give an ideal triangle in the above tessellation. We will then show how to
add bypasses to $\Sigma$ in such a way that the new dividing curves correspond to a different
triangle. After a sequence of such attachments we will eventually arrive at the triangle corresponding
to either $\{(1,1),$ $(0,1),$ $(1,2)\}$ or
$\{(1,0),(0,1),(1,1)\}$.

Before we begin the proof of Proposition \ref{normalization}, we need a preliminary lemma.

\begin{lem} The complete list of 
bypasses that can be attached to the initial configuration $\{v_1,v_2,v_1+v_2\}$ are as given in Figures
\ref{complist} and \ref{complist2} below.  The bypasses are subject to two conditions: (1) they do not give rise to destabilizations, and (2)
if $\delta\subset\Sigma_0$ is the attaching Legendrian arc for the bypass, and $p,q,r$ are intersections with
$\Gamma_0$ in consecutive order along $\delta$, then $p,q$ do not lie on the same dividing curve and $q,r$ do not
lie on the same dividing curve.
\end{lem}

The proof of the lemma is left to the reader --- one simply needs to check all the possibilities. To do this note
that the bypass can only involve two different dividing curves on $\Sigma$ (or else one gets a boundary-parallel
dividing curve on $\Sigma$). This observation will yield 12 possible dividing curve attachments, three of which 
will also result in a boundary-parallel dividing curve on $\Sigma.$ The nine remaining possibilities are listed in
Figures~\ref{complist} and \ref{complist2}.

\begin{figure}
	{\epsfxsize=3in\centerline{\epsfbox{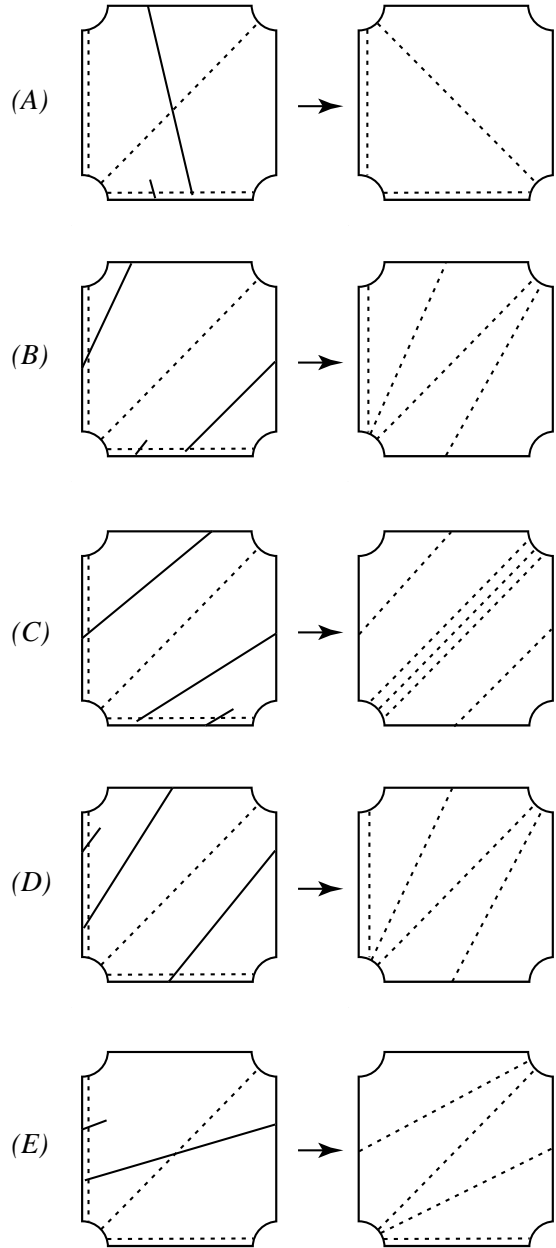}}}
	\caption{Allowable list of dividing curves.  The dotted lines are dividing curves, and 
	the Legendrian arcs of attachment for the bypasses are shown (solid lines). The bypasses are attached from
	the front. 
	The left side and right side of each tile are identified
	and the top and bottom are also identified.}
	\label{complist}
\end{figure}

\begin{figure}
	{\epsfxsize=3in\centerline{\epsfbox{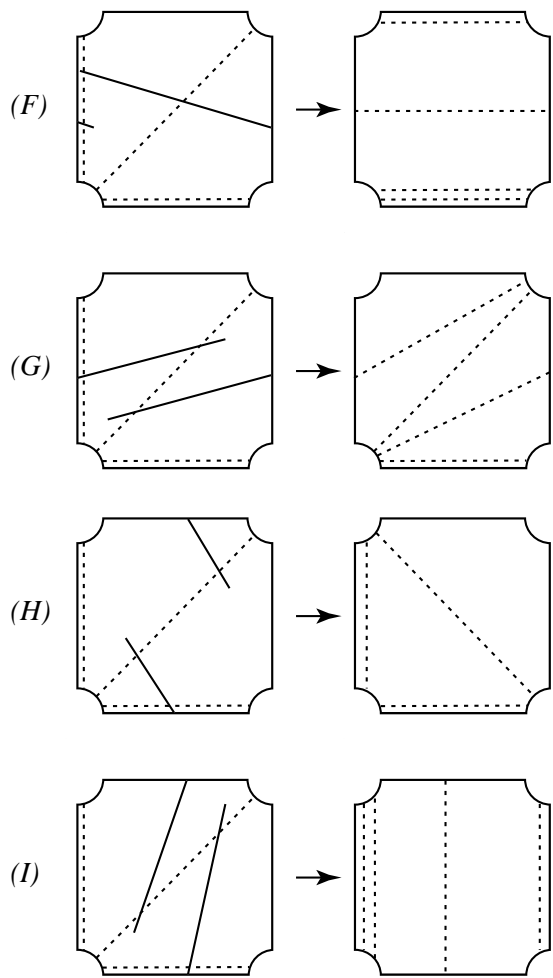}}}
	\caption{Continuation of allowable list of dividing curves.}
	\label{complist2}
\end{figure}
%

\begin{proof}[Proof of Proposition \ref{normalization}.] If $\Sigma_0$ has dividing set $\Gamma_0$,
then act via $\Psi$ to get $\Gamma_1=\Psi(\Gamma_0)$ on $\Sigma_1=\Psi(\Sigma_0)$,
which is a convex surface which is
isotopic to $\Sigma_0$.  Hence, we may freely modify
$\Gamma\mapsto \Psi^k(\Gamma)$ for any $k\in\Z$, via an isotopy.

Let $p, -{1\over p} \in \R\backslash \Q$ be the slopes corresponding to the two eigendirections of $\Psi$, with
$0< p< \infty$.  On the circle at infinity  $S^1_\infty$
of the standard tessellation of $D^2$, $p$ is an attracting fixed point and $-{1\over p}$ is a repelling fixed point
under the action of $\Psi$.  Suppose that the initial configuration $\{v_1,v_2,v_1+v_2\}$
corresponds to slopes $a,b,c$ respectively.  Then we may act repeatedly via $\Psi$ so
that $\infty\geq b>c>a \geq 0$, and we have three possibilities:
\be
\item $a,b,c>p$.
\item $a,b,c<p$.
\item $b,c>p$, $a<p$, or $b>p$, $a,c<p$.
\ee
Assume that in all the cases, bypass attachments do not give rise to boundary-parallel dividing curves (and hence
a destabilization), since that would contradict $\overline{\tb}=-3$.  This means that, for dividing curve configurations
of type III(a), we only use moves of type A through I in Figures \ref{complist} and \ref{complist2}.  For configurations of type
I(a), we only use the moves in Figure \ref{possibilities}.  To simplify notation, write $v_s$ for the shortest integral
vector with nonnegative entries, corresponding to a slope $0\leq s\in\Q$.

\s
\noindent
Case 1.  We may assume in addition that $1\leq a<c<b \leq \infty.$ Indeed, if $a\geq 1$ then we are done,
and if $a<1$ then $p<a,c,b<1,$ since any triangle in the tessellation whose vertices
are between $0$ and $\infty$ and whose clockwise-most vertex $<1$  has all its vertices $<1.$ Now
since  $\Psi(\infty)=1$ we may use $\Psi^{-1}$ to achieve the desired slopes.
We will show that using bypasses we can perform a sequence consisting of the following moves to our
$v_1, v_2,$ and $v_3=v_1+v_2$ (hence affecting $a,b,c$):
\begin{enumerate}
	\item replace $v_1,v_2,v_3$ with $v'_1, v'_2, v'_3$, where $v'_1=v_1-v_2, v'_2=v_2$ and $v'_3=v_1=v'_1+v'_2,$	
		(corresponds to right drawing in Figure \ref{move1})
	\item replace $v_1,v_2,v_3$ with $v'_1, v'_2, v'_3$, where $v'_1=v_1, v'_2=v_2-v_1$ and $v'_3=v_2=v'_1+v'_2,$
		(corresponds to left drawing in Figure \ref{move1}) or
	\item 
		replace $v_1,v_2,v_3$ with $v'_1, v'_2, v'_3$, where $v'_3=v_1,$ and $v'_1$ and $v'_2$ are
		determined by $v'_2=v'_1+v'_3$, and the triangle corresponding to $v'_1,v'_2,v'_3$ is the
		innermost
		triangle in the tessellation with $v'_3=v_1$ and all vertices counterclockwise of $1.$
		Here, `innermost' means closest to the center of the disk.
\end{enumerate}
Note that, if we start with a triangle $T$ corresponding to $v_1, v_2$ and $v_3=v_1+v_2$, then there is a
finite number $N_T$ of triangles that can be obtained from $T$ by the sequence of moves listed above. Moreover,
the triangle with vertices at $(1,2,\infty)$ is one of the triangles that can be so obtained. Now, if $T'$ is
obtained from $T$ by one of the moves above, then $N_{T'}<N_T.$ Thus any (sufficiently long) finite sequence of
these moves will take us from $T$ to the triangle with vertices at $(1,2,\infty).$

We now show how to find the bypasses that allow us to perform the above moves.
Since $\Psi(a)<1\leq a,b,c$, $|\Psi(v_a)\cap (v_a\cup v_b\cup v_c)|>2$, and, by using an annulus
$d\times I$ with $d$ parallel to $\Psi(v_a)$, we find a bypass of type A, F, or H in Figures
\ref{complist} and \ref{complist2}.  The bypasses are of the stated type, since the slope of
the bypass is smaller than $a,b,$ or $c.$ Suppose we do not already have
$\{a,b,c\}=\{1,2,\infty\}$.  Moves of type A and H replace the ideal triangle with
vertices $a,b,c$ by another triangle with vertices $a,b,c'$ (see Figure \ref{move1}).  Here $c'$ corresponds
to the slope of $\pm(v_2-v_1)$.  Thus, $v_1, v_2$ and $v_3$ are affected by applying one of the first
two moves mentioned above.
\begin{figure}
	{\epsfysize=1.5in\centerline{\epsfbox{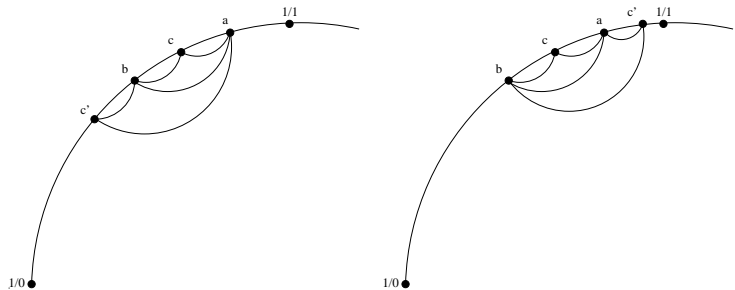}}}
	\caption{Type A and H moves.}
	\label{move1}
\end{figure}

A move of type F will give rise to a configuration of type I(a) with $n=3$ and $m=1$.  Just
as we denoted a configuration of type III(a) by $\{a,b,c\}$, we will denote a configuration
of type I(a) by $\{a\}$.   We deal with the case $a=1$ at the
end of the paragraph.  Otherwise, in order to modify
$\Gamma_0= \{a\}$, take an annulus $d\times I$, where $d$ has slope $\Psi(a)$.
Note that $p< \Psi(a)<1$.
Since bypass attachments cannot give rise to boundary-parallel dividing curves on $\Sigma$,
there are two possibilities, given in Figure~\ref{possibilities}.  Here $b'$ is the smallest
slope for which there exists an edge of the tessellation from $a$ to $b'$. (Note: in general there are
other possibilities for the slope of the bypass but with this choice of basis for the torus and our choice
of $d$ only those shown in Figure~\ref{possibilities} are possible.)
\begin{figure}
	{\epsfxsize=1.5in\centerline{\epsfbox{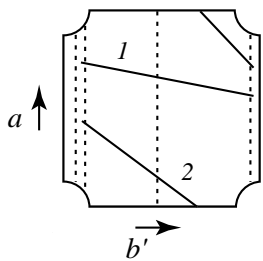}}}
	\caption{Two possibilities for bypass attachment.}
	\label{possibilities}
\end{figure}
After the bypass attachment, $\Sigma$ has dividing curves of type III(a) and two of the dividing curves
have slopes $a, b'$ respectively. Let $v_1$ and $v_2$ be the minimal integral vectors corresponding
to these slopes. Then the third dividing curve on $\Sigma$ is spanned by $v_3=v_1+v_2$ or $v_2-v_1.$
In the first case we have performed a move of type 3 listed above and in the second case we have
performed a move of type 3 followed by one of type 2.
Thus we may perform a sequence of the above moves eventually 
obtain $\{1,2,\infty\}$ or at some point we obtain a configuration of type I(a) with curves of
slope $1$.  
In this case use $\Psi^{-1}$ to modify $\{1\}\mapsto \{\infty\}$.  Let this be
$\Gamma_0$.  Then, using the annulus $d\times I$ where $d$ has slope $\Psi(\infty)=1$,
we obtain a bypass which modifies $\{\infty\}\mapsto \{1,2,\infty\}$.

\s
\noindent
Case 2.  Assume $a<c<b<p$.  Without loss of generality, assume
$0\leq a<c<b \leq{1\over 2}$ (this argument is the same as the one at the beginning of 
Case 1).
Use the annulus $d\times I$, where $d$ has slope
$\Psi(b)$, to obtain a bypass of type A, F, or H.  The same argument as in Case 1 allows us
to isotop $\Sigma$ so that the dividing curves are of type I(a) with slope 0 or of type
III(a) with slopes $0,{1\over 2},{1\over 3}.$ In the first case, attaching the only possible bypass along
$\Psi(0)={1\over 2}$ which does not yield a
boundary-parallel component on $\Sigma_0$, we obtain $\{0\}\mapsto \{0,1,\infty\}$.  
Next suppose we have $\{0,{1\over 2},{1\over 3}\}$.  If we use the annulus $d\times I$ where $d$ has
slope ${1\over 2}$, then $\Psi({1\over 2})={3\over 5}$, and F cannot happen (for F to happen, we need
$\Psi({1\over 2})>1$).  Therefore, $\{0,{1\over 2},{1\over 3}\}\mapsto \{0,{1\over 2}, 1\}$, which is 
dealt with in Case 3.

\s
\noindent
Case 3.  Here we start with a triangle $T$ in the Farey tessellation whose vertices straddle $p.$
Below we show that attaching bypasses will either put us back in Case 1, which we have already
dealt with, or we perform one of the following moves on $v_1, v_2$ and $v_3=v_1+v_2:$
\begin{enumerate}
	\item replace $v_1,v_2,v_3$ with $v'_1=v_1, v'_2=v_1+v_2$  and $v'_3=v_2+2v_1$ (if $a<p$ and $b>c>p$) or
	\item replace $v_1,v_2,v_3$ with $v'_1=v_1+v_2, v'_2=v_2$ and $v'_3=v_1+2v_2$ (if $a<c<p$ and $b>p$).
\end{enumerate}
Assuming this is true for the moment. There is some $k$ such that the vertices 
$\{\Psi^k(0), \Psi^k(1), \Psi^k(\infty)\}$ span a triangle $G$ in the Farey tessellation
that is disjoint from $T$ and separates $T$ from $p.$  Each time we perform a move
of type 1 or 2 to $T$, we obtain a new triangle $T'$ lying between $T$ and $G$.
Since there are only a finite number
of triangles in the tessellation lying between $T$ and $G$, this process must stop. It can only
stop if $T=G$ (or a bypass attachment puts us in Case 1), in which case we are done.

Suppose that $b>c>p$ and $a<p$.  Use the annulus $d\times I$, where
$d$ has slope $\Psi(c)$, to obtain a bypass of type C, E, or G.  C gives
$\{a,b,c\}\mapsto \{c\}$, with $c>p$, which was already treated under Case 1.  
Bypasses of type E or G modify $v_1,v_2,v_1+v_2$ to $v_1,v_1+v_2,2v_1+v_2$ --- and we are still in
Case 3 (though we might have $b>p$ and $a<c<p$). 
Now suppose that $b>p$ and $a<c<p$.  Similarly, we have bypasses of type B, D, or I.  
I puts us in Case 1, and B, D modify $v_1,v_2,v_1+v_2$ to $v_1+v_2,v_2,v_1+2v_2$ (still in 
Case 3).
\end{proof}

\begin{prop} \label{g2} (1) There exists no tight contact structure on $\Sigma\times I$ for which
$\Gamma_0$ is of type III(a) with slopes $\{0,1,\infty\}$.  (2) There exists exactly
one tight contact structure on $\Sigma\times I$ for which $\Gamma_0$ is
of type III(a) with slopes $\{1,2,\infty\}$, and $(\Sigma_0)_\pm$
 --- the signs of each region divided by $\Gamma_0$ --- are fixed.
\end{prop}

\begin{proof}
Note that $H=\Sigma\times I$ is a genus two
handlebody.  In general, when analyzing handlebodies with convex boundary, we use compressing disks
$D_1$ and $D_2$ so that the cut-open manifold is a 3-ball.  To obtain a smooth convex boundary of 
$H=\Sigma\times I$, we need
to round the edges $\bdry\Sigma \times\{0,1\}$.  This is done using the Edge-Rounding Lemma (Lemma~\ref{edge}).
When using Edge-Rounding, we must be careful to remember that, since $\bdry\Sigma\times[0,1]$ was the convex
boundary of a neighborhood of a Legendrian curve, there is holonomy (as we go from $\Sigma_0$ to $\Sigma_1$).   
We then arrange for $D_1$, $D_2$ to have Legendrian boundary --- for this we use the Legendrian Realization
Principle.  Then, $D_1$, $D_2$ with Legendrian boundary are perturbed, fixing the boundary, so they become convex.
What remains is to study the configuration of dividing curves on each compressing disk.  

We have $\partial(\Sigma\times [0,1])=\Sigma_0\cup\Sigma_1\cup A$ where $A$ is the annulus 
$\partial \Sigma\times [0,1].$ Below we discuss the dividing curves on $\Sigma_0$ and $\Sigma_1$ but first we 
consider the dividing curves on then annulus $A$ and how they are related to the dividing curves on the
$\Sigma_i$'s. To this end consider Figure~\ref{fig:monodromy}. In this figure we
identify the right and left edges, then the center rectangle forms the annulus $A.$ The 
top and bottom parts of the picture form neighborhoods of $\partial\Sigma_i$ in $\Sigma_i$ 
and the shaded region is obtained using the Edge-Rounding Lemma. The dividing curves on $\Sigma_i$
divide the $\partial\Sigma_i$ into 6 intervals. We label the intervals 1-6 so we can identify them
with the corresponding intervals on the annulus $A.$ 
\begin{figure}
	{\epsfysize=2.3in\centerline{\epsfbox{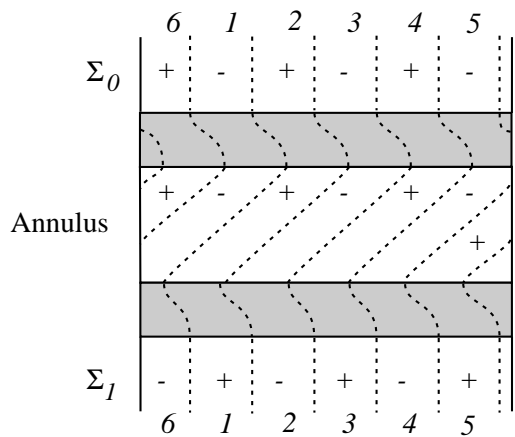}}}
	\caption{Dividing curves (dashed lines) on, and near, the annulus $A.$}
	\label{fig:monodromy}
\end{figure}

(1) $\Gamma_0= \{0,1,\infty\}$ and $\Gamma_1= \{ {1\over 2},{2\over 3},1\}$.
See Figure~\ref{notight}.  There exists a disk $D=\delta\times[0,1],$ where $\delta$ is an arc with slope 1, with 
geometric intersection $|\bdry D\cap \Gamma_{\bdry H}|=0$. Note that $\delta\times\{0,1\}\subset
\partial D$ are represented by solid lines in Figure~\ref{notight}, and
$(\bdry \delta)\times I\subset \partial D$
lies on $A=\partial \Sigma\times[0,1].$ Though $(\bdry\delta) \times I$ is not drawn in Figure~\ref{fig:monodromy},
one may easily determine that $\partial D\cap \Gamma_{\bdry H}=\emptyset$.
By the Legendrian Realization Principle, we may realize $\bdry D$ as a zero-twisting (rel $\bdry H$)
Legendrian curve.  This implies that there is an overtwisted disk and that
there cannot be a tight contact structure on $H$ with this boundary
condition.

\begin{figure}
	{\epsfxsize=4.5in\centerline{\epsfbox{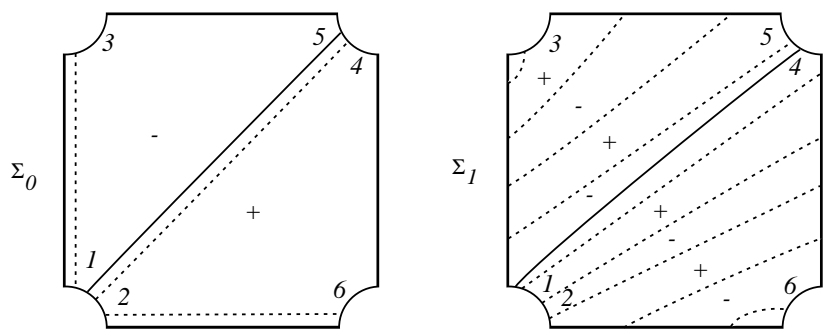}}}
	\caption{Dividing curves on $\Sigma_0$ and $\Sigma_1$ (dashed lines). The solid line 
	is our choice of $\delta.$}
	\label{notight}
\end{figure}


(2)  $\Gamma_0= \{1,2,\infty\}$ and $\Gamma_1= \{ {2\over 3},{3\over 4} ,1\}$.
Let $\delta_1$ be an arc on $\Sigma$ with slope 1 and let $\delta_2$ be an arc with infinite slope.
Choose the compressing disks $D_1$ and $D_2$ for $H$ 
to be isotopic to $\delta_1\times I$ and $\delta_2\times I$, drawn as in Figure~\ref{yestight}.
After using the Legendrian Realization Principle, we take the boundaries 
to be Legendrian and $\partial D_1$ and $\partial D_2$ to intersect the dividing
curves on $\partial H$ (minimal geometric intersection number) 2 and 6 times, respectively.

\begin{figure}
	{\epsfxsize=4.5in\centerline{\epsfbox{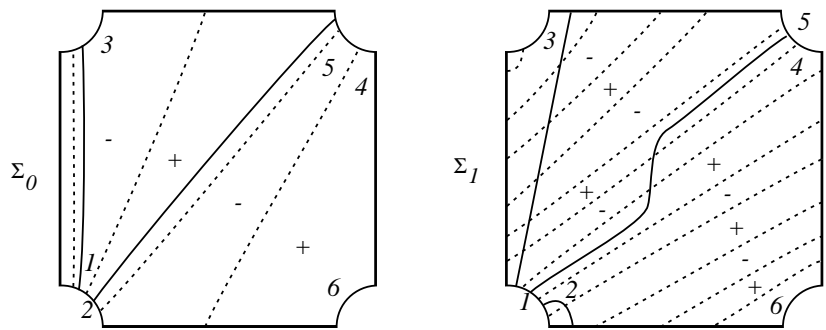}}}
	\caption{Dividing curves on $\Sigma_0$ and $\Sigma_1$ (dashed lines). The solid line 
	is our choice of $\delta_i.$}
	\label{yestight}
\end{figure}

There is only one possibility for the dividing curves on $D_1$.
However, there are several possible dividing curve
configurations for $D_2$.
Consider the intersections between $\delta_2\times\{1\}$ and $\Gamma_1$ --- count the intersections
from bottom to top along $\delta_2\times\{1\}\subset \Sigma_1$ shown in Figure~\ref{yestight}.
One may easily check that there cannot exist a bypass
along $\partial D_2$ which straddles the second intersection or the fourth intersection,
without immediately yielding a bypass for $K$. Thus, the only two possibilities for $D_2$ are shown in
Figure~\ref{gandb}. Now, if $D_2$ had the right-hand configuration shown in Figure~\ref{gandb}, then a
bypass straddles the third intersection between $\delta_2\times\{1\}$
and $\Gamma_1$. Note this
bypass is nested inside another bypass on $D_2.$  If we added
both these bypasses in succession, the resulting copy of
$\Sigma$ would have dividing curves $\{0,1,\infty\}.$ Hence
we could find an overtwisted disk as above. Thus, the dividing curves on $D_2$ are shown on the left-hand side
of Figure~\ref{gandb}.

\begin{figure}
	{\epsfxsize=3.5in\centerline{\epsfbox{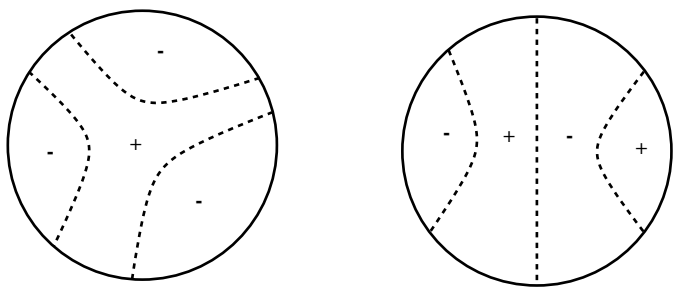}}}
	\caption{Possible dividing curves on $\Delta_2.$}
	\label{gandb}
\end{figure}

Now, since we have normalized each of $D_1$, $D_2$ to have a unique dividing curve configuration, 
we have a unique tight contact structure on $H$ with given configurations $\Gamma_0$, $\Gamma_1$,
up to isotopy.  This is because any two tight contact structures with these boundary conditions can
first be matched up along $D_1$, $D_2$, and then matched up inside $H\backslash (D_1\cup D_2)$ using 
Theorem~\ref{thm:uniqueB3}.
\end{proof}

We have normalized $\Gamma$ on $\Sigma$ to be of type III(a) with slopes $\{1,2,\infty\}$.
To finish the proof of Theorem~\ref{thm:uniquemax}
we just note that we can also normalize the signs of the two regions of $\Sigma\backslash \Gamma$ by
acting via $-id$ on $\Sigma$.  Since $-id\circ \Psi\circ (-id)^{-1}=\Psi$, we may use this 
{\it isomorphism} to switch $\Sigma_+$ and $\Sigma_-$.  Therefore, we obtain a contactomorphism 
$f:S^3\rightarrow S^3$ taking $K$ to $K'$ as discussed after the statement of Theorem~\ref{thm:uniquemax}.

\subsection{Destabilization}

\begin{prop}\label{prop:stabilize}
	If $K$ is an oriented Legendrian figure eight knot and $\tb(K)=-3-k<-3,$
	then there exist positive integers $k_1$ and $k_2$ such that $k=k_1+k_2,$
	$r(K)=k_2-k_1$ and
	$K=S_-^{k_1}(S_+^{k_2}(K'))$, where $K'$ is the Legendrian figure eight knot with
	maximal Thurston-Bennequin invariant.
\end{prop}

\proof
Suppose $K$ is a Legendrian figure eight knot with $\tb(K)<-3.$ From Lemma~\ref{lem:invts}
we know that there must be more than 3 dividing curves on $\Sigma.$ Thus Proposition~\ref{prop:curves}
implies that the dividing curves on the Seifert surface $\Sigma$ can be normalized so that we have
I(a) (with $n\geq 5$), II(a), or III(a). In each of the cases, we prove that either we can normalize 
$\Gamma$ into a standard form, or there exists a destabilization --- this is done in the same way as in
Proposition \ref{normalization}.  Then, in the same way as in Proposition \ref{g2}, we prove that there always
exists a destabilization, when $\Gamma$ is in the standard form.  

\begin{figure}
	{\epsfxsize3in\centerline{\epsfbox{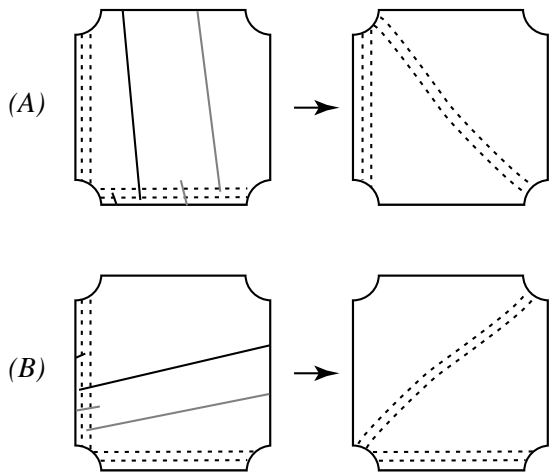}}}
	\caption{Allowable bypasses. The Legendrian arcs of attachment are shown.  The bypasses are attached
	along solid black or grey lines from the front.}
	\label{fig9}
\end{figure}

\s
\noindent
Case II(a).  Refer again to the tessellation picture, and denote $\Gamma_0$ by $\{a,b\}$,
where $\infty\geq b>a\geq 0$ (without loss of generality) are the slopes of the two isotopy classes of arcs on
$\Sigma_0$.  Figure \ref{fig9} gives a complete list of bypasses 
that can be attached, without immediately giving rise to a destabilization.  
First assume $\infty\geq b>a\geq 1$, and that $\{a,b\}\not = \{1,\infty\}$.
If $d$ is a closed curve with slope $a$, then use $\Psi(d)\times I$ to
obtain a bypass along $\Sigma_0$, using the Imbalance Principle.  If there is no destabilization, then 
the bypass attachment is of type A, and 
$\{a,b\}\mapsto \{a',b\}$, where $a>a'\geq 1$.
If we repeat this procedure, we either arrive at $\{1,\infty\}$ or obtain a destabilization if $a$ is already the smallest
slope for which there exists an edge of the tessellation to $b$.
Finally, if we take a closed curve $d$ with slope $1$, then $\Psi(d)$ has slope ${2\over 3}$.  The bypass
along $\Sigma_0$ with $\Gamma_0= \{1,\infty\}$, which comes from $\Psi(d)\times I$, cannot be of
type A or B, and hence we have a destabilization.

Next assume $\infty\geq b>p$ and $p>a\geq 0$.  Define $c$ to be the slope corresponding to $v_a+v_b$.  If $c<p$,
then use $\Psi(d)\times I$, where $d$ has slope $b$.  The resulting bypass is not of type A or B, 
since the slope of $\Psi(d)>c>0$ while the bypasses of type A and B have slope less than $c.$ Thus
we can destabilize.  If $c>p$, then use $\Psi(d)\times I$, where $d$ has slope $a$.
A type B move modifies $\{a,b\}\mapsto \{a,c\}$.  We keep repeating, until eventually
the new $c$ satisfies $c<p$.

If $p>b>a\geq 0$, then we may assume ${1\over 2}\geq b>a\geq 0$.  Assume $\{a,b\}\not =\{0,{1\over 2}\}$
already.    The situation is similar to $\infty\geq b>a\geq 1$
above, and we eventually obtain $\{0,{1\over 2}\}$ or a destabilization.
A bypass attachment of type A (which must exist if there
is no destabilization) will put us into $\{0,1\}$ which is already done.
Therefore, Case II(a) always destabilizes.

\begin{figure}
	{\epsfxsize3in\centerline{\epsfbox{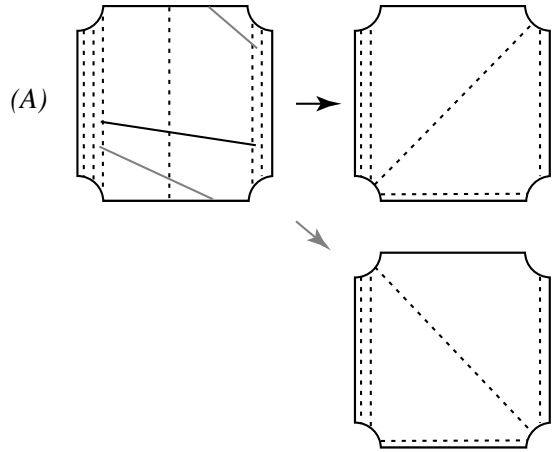}}}
	\caption{Allowable bypasses.  Case I(a). }
	\label{fig10}
\end{figure}

\begin{figure}
	{\epsfxsize3in\centerline{\epsfbox{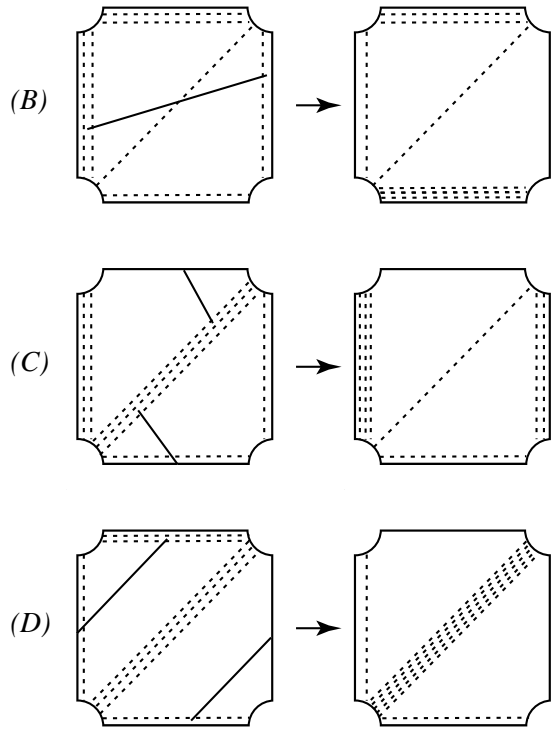}}}
	\caption{Allowable bypasses. Case III(a), and two of the isotopy classes have more than one parallel copy. }
	\label{fig11}
\end{figure}

\begin{figure}
	{\epsfxsize3in\centerline{\epsfbox{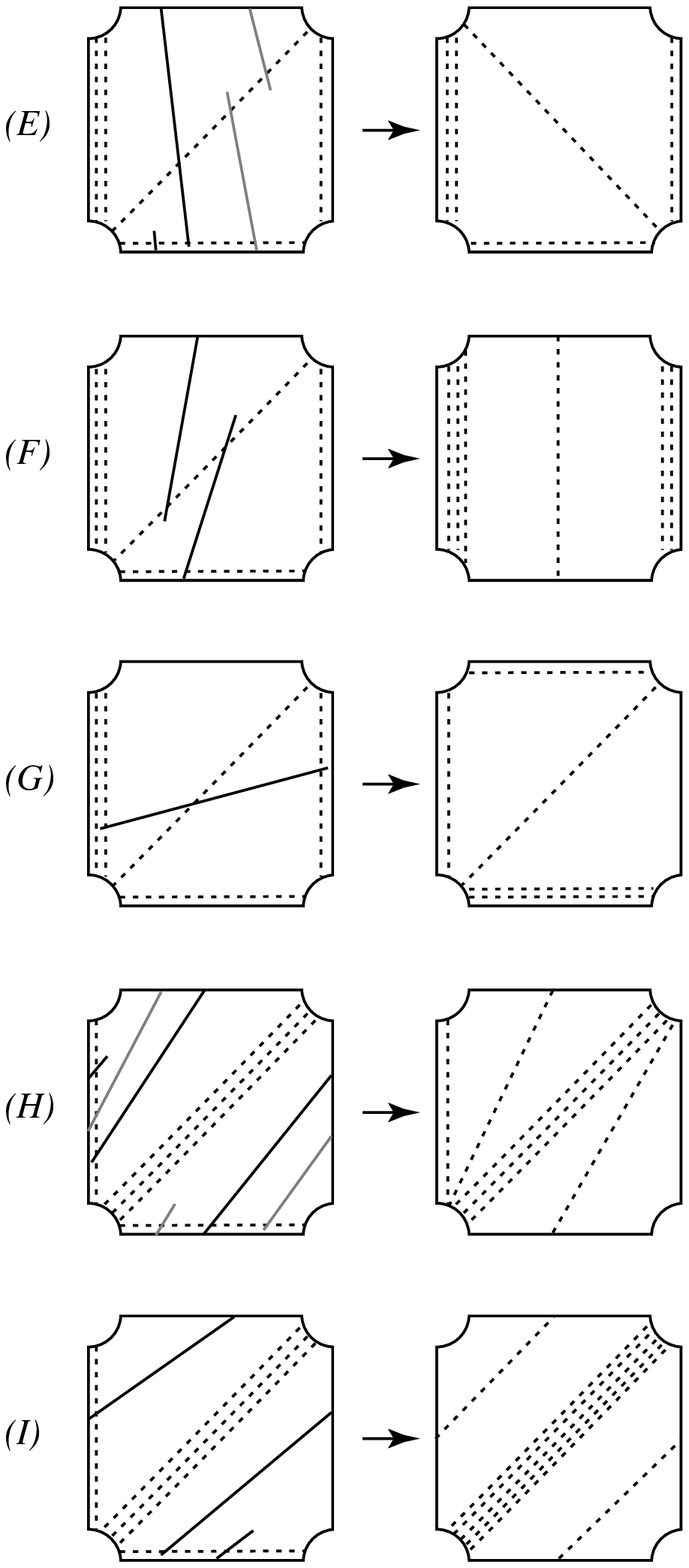}}}
	\caption{Allowable bypasses. Case III(a), and two isotopy classes have one copy each. Bypasses are attached
	along solid black or grey lines.}
	\label{fig12}
\end{figure}

\begin{figure}
	{\epsfxsize3in\centerline{\epsfbox{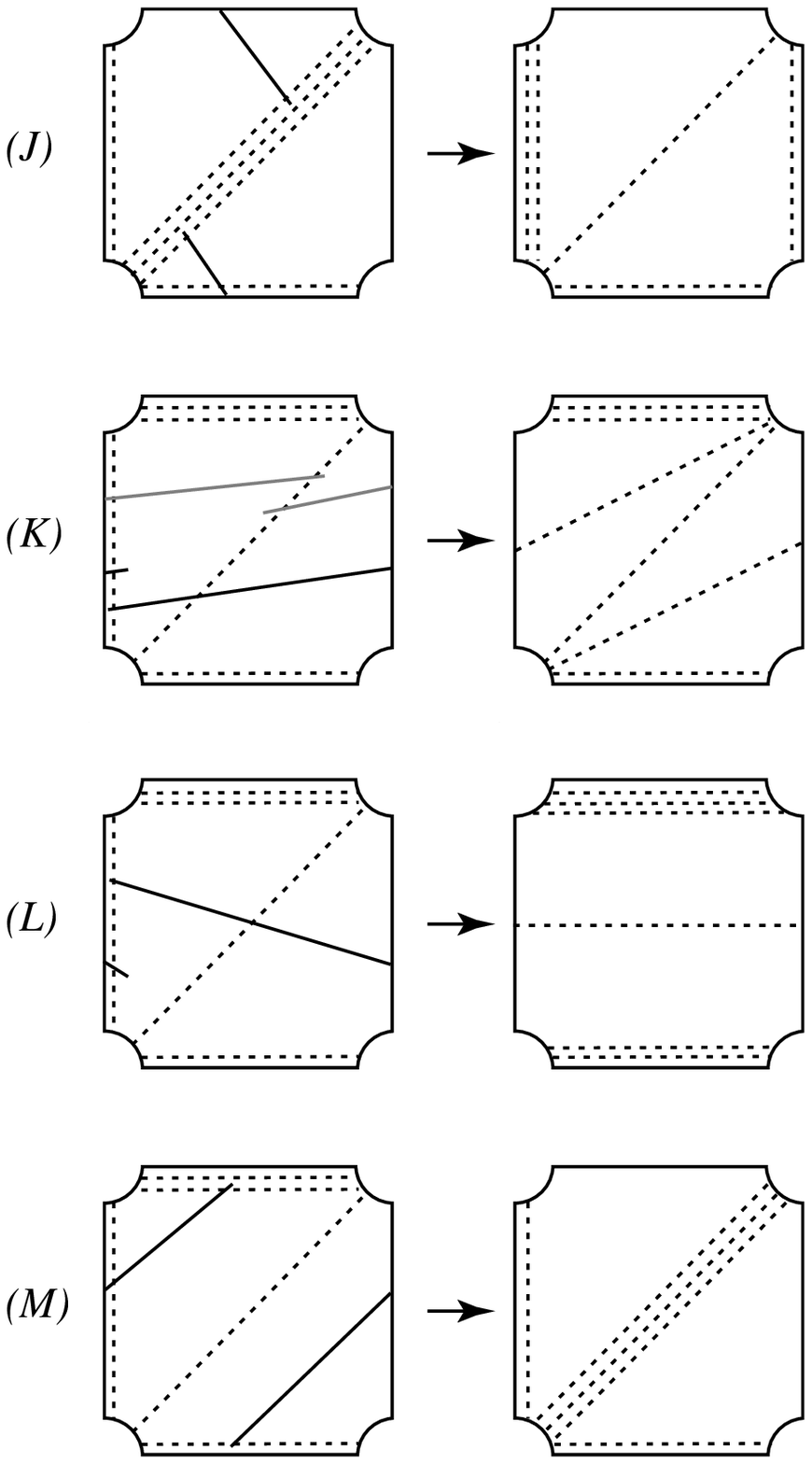}}}
	\caption{Allowable bypasses. Case III(a), and two isotopy classes have one copy each. Bypasses are attached
	along solid black or grey lines.}
	\label{fig122}
\end{figure}

\s
\noindent
Case III(a).  Consider Figures \ref{fig10}, \ref{fig11}, \ref{fig12} and \ref{fig122}, which list all the possible
bypasses.

Figure~\ref{fig10} (just Type A) is for I(a) (Note: as discussed in the proof of Proposition~\ref{normalization}
there are other possibilities for the slopes of the bypasses in I(a) but with the proper choice of basis
for the torus we can always arrange to be in the situation shown in Figure~\ref{fig10}).
When two of the isotopy classes of arcs have more than one parallel copy, then we have Figure \ref{fig11} and
one of
B, C, or D.  When two isotopy classes have one copy each, then we have Figures \ref{fig12} and \ref{fig122} and one of
E, F, G, H, I, J, K, L, or M.  As before, represent a I(a) configuration by $\{a\}$ and a III(a) configuration by $\{a,b,c\}$.  Let
$\m(s)$ be the multiplicity (number of parallel copies) of the isotopy class with slope $s$.

Assume we have III(a) and the Legendrian curve does not destabilize.
Notice that we may assume that two of the isotopy classes will have $mult=1$ after
repeated attachments of bypasses  --- this follows from the fact that
if $mult>1$ for 2 isotopy classes, then Figure \ref{fig11} implies that there are only moves B, C, and D, and these
yield two isotopy classes of $mult=1$ after repeated application (note depending on the situation one can only
find bypasses of type B, C or D and not some combination of these).
Therefore we assume that two classes have $mult=1$.

Now assume, in addition, that $\infty \geq b>c>a\geq 1$, and
$\{a,b,c\}$ is not already $\{1,2,\infty\}$.  As before, take a closed curve $d$ with slope $a$,
and consider $\Psi(d)\times I$.  Since $\Psi(d)$ intersects each curve of $\Gamma_0$ at least once, there exists
a bypass of slope  $\Psi(a)$ with $1> \Psi(a)>p$.  The only possibilities are then E, G, J, and L.
We also want to distinguish between $\Psi(a)>a'$ and $\Psi(a)<a'$, where $a'$ corresponds to $v_a-v_b$.
If $\Psi(a)>a'$, then we have G, J, or L.  If $\Psi(a)<a'$, then we have E or J.

Assume $\Psi(a)<a'$.
Suppose first $mult(a)=mult(b)=1$.  Then we have move J, and repeated application yields
$mult(a)=mult(c)=1$. Using move E, we obtain $\{a,b,c\}\mapsto \{a',a,b\}$ ($\infty\geq b>a>a'\geq 1$),
$mult(a')=mult(a)=1$ and $mult(b)>1$.  $\{a',a,b\}$ corresponds to an ideal triangle which
lies further inward on $D^2$  --- the kind on the
right-hand side of Figure \ref{move1}.

Assume $\Psi(a)>a'$.
Suppose first $mult(a)=mult(b)=1$.  Then repeated application of move J yields $mult(a)=mult(c)=1$.
Next, repeated application of move G yields $mult(c)=mult(b)=1$.  Then we are now left with
L which maps $\{a,b,c\}\mapsto
\{a\}$.  Taking $d$ with slope $a$ and annulus $\Psi(d)\times I$,
we get $\{a\}\mapsto \{a'',a,b''\}$, where $a''\geq 1$ is the smallest slope with an edge to
$a$, $mult(a'')=mult(b'')=1$, and either $a''<b''<a$ or $a''<a<b''$.
Thus, eventually the smallest slope $a$ of $\{a,b,c\}$ goes to 1 and  we get to $\{1,2,\infty\}$
with $mult(1)=mult(2)=1$, $mult(\infty)
>1$, or there exists a destabilization.    Once we have $\{1,2,\infty\}$, use $\Psi(d)\times I$, where $d$
has slope $1$, and repeatedly apply J and G until $mult(1)>1$, $mult(2)=mult(\infty)=1$, followed by L, which gives
us $\{1,2,\infty\}\mapsto \{1\}$.  Since we may map via $\Psi^{-1}$ to get $\{\infty\}$, we may use $\{\infty\}$
or $\{1,2,\infty\}$ as we wish.

\begin{figure}
	{\epsfysize1in\centerline{\epsfbox{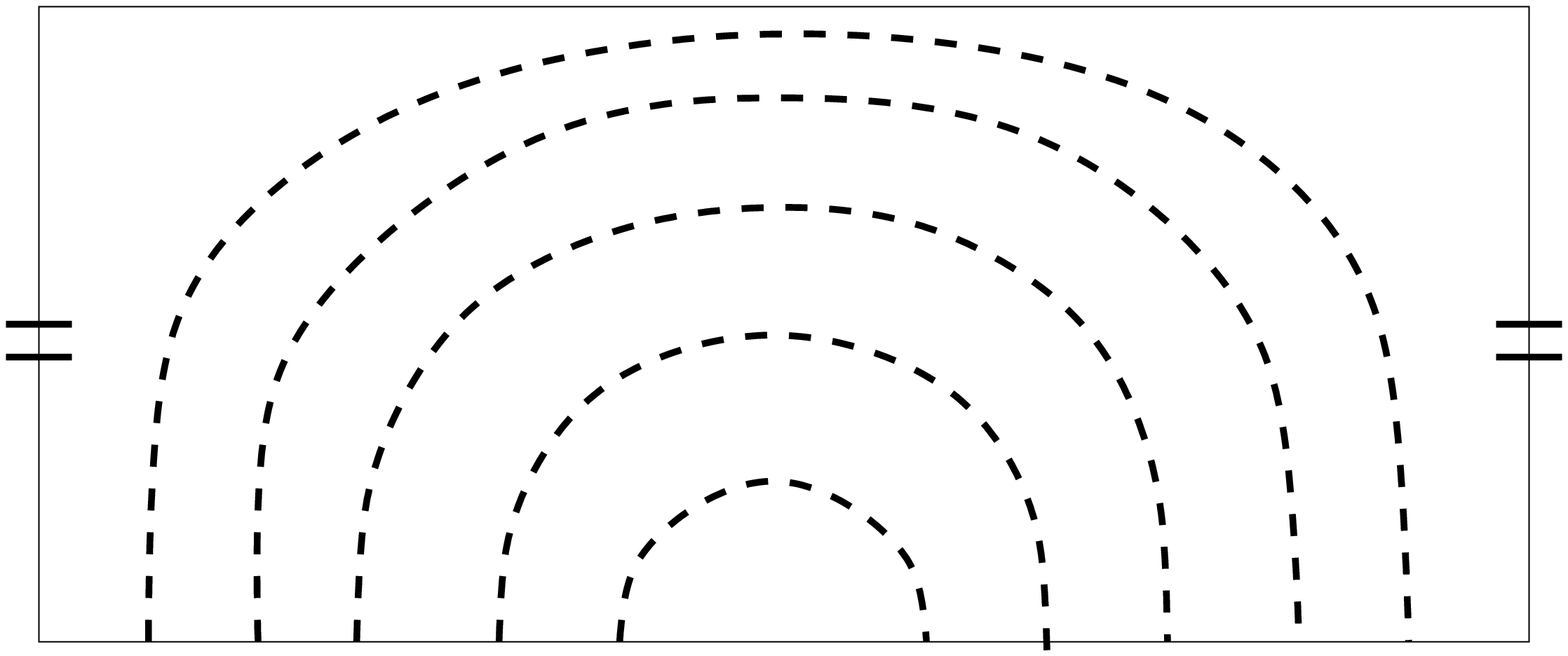}}}
	\caption{Nesting.  The sides of the annulus are identified.  The top is $d\times\{1\}$ and the
	bottom is $d\times \{0\}$. }
	\label{fig13}
\end{figure}

\begin{figure}
	{\epsfxsize4.5in\centerline{\epsfbox{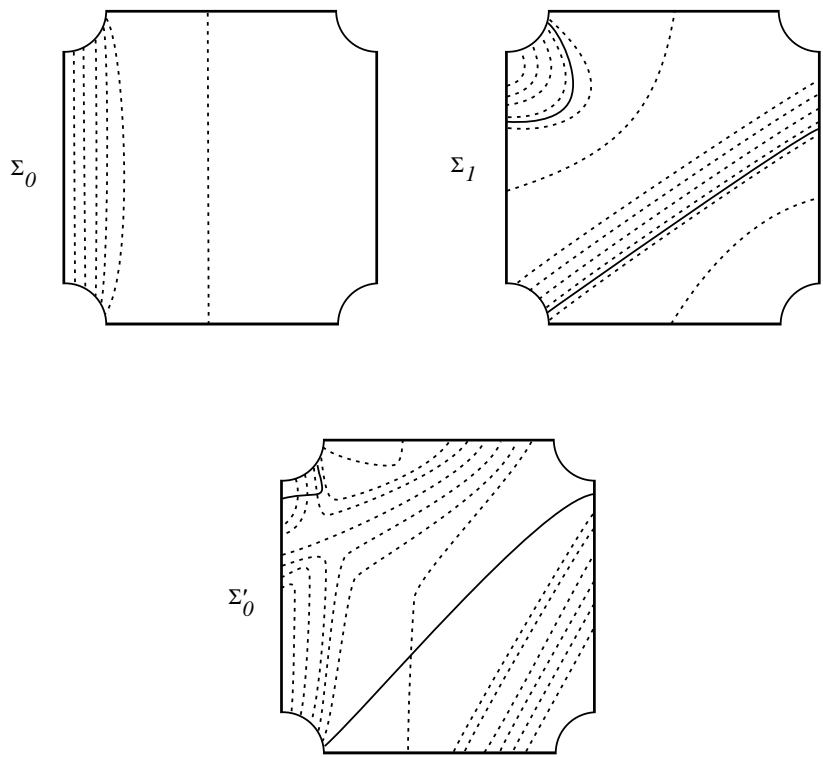}}}
	\caption{$\Sigma_0, \Sigma'_0$ and $\Sigma_1.$ The non-dashed curve is $\delta.$}
	\label{fig14}
\end{figure}

Up to now we made no use of the holonomy on $\bdry \Sigma\times I$.  We now mimic the
proof of Proposition \ref{g2} and study the genus two handlebody $H=\Sigma \times I$.   $\Sigma_0$ has
$\Gamma_0$ consisting of $n$ (odd) arcs and $1$ closed curve, all of slope $\infty$.
$\Sigma_1$ has $\Gamma_1$ consisting of $n$ arcs and $1$ closed curve, all of slope $1$.  Take
$d\times I$, where $d$ has slope $1$.  The Legendrian realization $d\times\{1\}$ has twist number 0 on $\Sigma_1$,
whereas the Legendrian realization $d\times\{0\}$ has twist number $n+1$ on $\Sigma_0$.  Since there is only
one allowable bypass attachment onto $\Sigma_0$ with slope $1$, all the dividing curves on $d\times I$ must
be nested `in parallel'.   See Figure \ref{fig13}.
If we attach all the nested bypasses in succession to $\Sigma_0$, we eventually obtain $\Sigma_0'$, which is
identical to $\Sigma_1$, except for spiraling.  See Figure \ref{fig14}.
If we take the meridional disk $D=\delta\times I$,  where $\delta$ is an arc with slope $1,$ then
$|\bdry D\cap \Gamma_{\bdry H}|=n-1$.  As long as $n>3$, there will exist a bypass along $\bdry D.$
In fact, since $D$ is a disk, at least one bypass will lie entirely on $\Sigma_0.$ Since
the bypass does not involve the closed dividing curve on $\Sigma'_0$, it will have to lead to a boundary-parallel
dividing curve on $\Sigma'_0.$ Thus we have a destabilization of our knot.

Next assume ${1\over 2}\geq  b>c>a\geq 0$.  As above, there is a destabilization or we get to
$\{0,{1\over 3}, {1\over 2}\}$, with $mult(0)=mult({1\over 3})=1$, $mult({1\over 2})>1$, or to
$\{0\}$.  In the former case, a bypass of
type E then gives $\{0,{1\over 2},1\}$.  In the latter case, a type A move gives $\{0,{1\over 2},1\}$.

Assume $p>a\geq 0$, $\infty \geq  b>p$.  First assume $b>c>p$.   Then we let $d$ have slope $a$, and consider
$\Psi(d)\times I$.  Here $\Psi(a)<p$.  Then we have G, I, K, or M.    I gives $\{c\}$ with $c>p$, which is done.
If we have M, then the next bypass operation must be I, and we are done.  If we have G, then the next step is
M, and we are done again.    We are left with K.   Next assume $p>c>a$.  Then let $d$ have slope $b$, and consider
$\Psi(d)\times I$.  Here $\Psi(b)>p$.  We will have F, H, J, or M.    As before, F gives $\{b\}$ with $b>p$, which is done.
J must be followed by F, and M must be followed by H or J, so we are left with H.
At the end of the day, we are left with $\{\Psi^k(0), \Psi^k(1),\Psi^k(\infty)\}$ (or simply $\{0,1,\infty\}$),
and $mult(0)>1$, $mult(1)=mult(\infty)=1$.  Now (see Proposition \ref{g2}), a handlebody $H$ cut open along
$\Sigma$ with this configuration must give rise to a destabilization --- use the disk $D=\delta\times I$, where
$\delta$ has slope $\infty$.

\s
\noindent
Case I(a). Given the configuration $\{a'\}$ with $\infty>a'>1$, move A will put us into case III(a) with
$\infty\geq b>c>a\geq 1$, which is done.  Given $\{a'\}$ with ${1\over 2}>a'>0$, move A will put us into case
III(a) with ${1\over 2}\geq b>c>a\geq 0$, which is also done.
\qed

\s
\noindent
We have now proved that for $\mathcal{K}$, the isotopy class of the figure eight knot, there exists
a unique maximal Thurston-Bennequin representative, and any Legendrian figure eight knot destabilizes to this
unique representative.  Now, using Lemma \ref{aboutstabilization}, we find that the classification easily extends
to non-maximal $\tb$ Legendrian knots.  This proves Theorem \ref{thm:maineight}.      \qed

\subsection{Transversal Figure Eight Knots}
Theorem \ref{thm:maineight} implies that the figure eight knot is stably simple and thus we may
use Theorem~\ref{thm:stably-transversally} to conclude
\begin{cor}
	If $K$ and $K'$ are two transversal figure eight knots then they are transversally  isotopic
	if and only if $l(k)=l(K').$ Moreover, the self linking numbers of transversal figure eight knots
	realize precisely the set of odd integers $\leq -3.$
\end{cor}

\s\s
\noindent
{\em Acknowledgments:}
The first author gratefully acknowledges the support of an NSF 
Post-Doctoral Fellowship(DMS-9705949) and Stanford University.    The second author would
like to thank Will Kazez and Gordana Mati\'c for a stimulating research environment at the University of
Georgia.


\end{document}